\documentclass{article}
\usepackage{graphicx, subfigure, amsfonts, amsmath, amsthm, tikz, url}
\usepackage{a4, color}
\usepackage[top=1in, bottom=1.25in, left=1.25in, right=1.25in]{geometry}
\newcommand{\ignore}[1]{}
\newcommand{\paf}[2]{\frac{\partial #1}{\partial #2}}
\newcommand{\naf}[2]{\frac{\text d #1}{\text d #2}}

\newcommand{\discretization}{discretization}
\newcommand{\discretizations}{discretizations}
\newcommand{\adjointness}{adjointness}

\newcommand{\diag}{\mbox{diag}}

\newcommand{\eps}{\varepsilon}

\newcommand{\correction}[2]{{ \color{myorange} $<<${#1}\footnote{#2}$>>$ }}

\newcommand{\corr}[2]{{ \color{myorange} $<<${#1}$>>$}{\cororig{#2} }}
\renewcommand{\correction}[2]{#1} \renewcommand{\corr}[2]{#1}
\begin{document}
\title{Symmetry-preserving finite-difference discretizations of arbitrary order on structured curvilinear staggered grids}
\author{Bas van 't Hof\footnote{Corresponding author. bas.vanthof@vortech.nl. VORtech, Westlandseweg 40d, 2624AD Delft.} \qquad Mathea~J. Vuik\footnote{thea.vuik@vortech.nl.}}

\maketitle

\begin{abstract}
Symmetry-preserving (mimetic) discretization aims to preserve certain properties of a continuous differential operator in its discrete counterpart.  
  For these discretizations, stability and (discrete) conservation of mass, momentum and energy are proven in the same way as for the original continuous model.  
  
  This paper presents a new finite-difference symmetry-preserving
  \correction{space}{Reviewer 1, Remark 3: semidiscrete, shocks} discretization. \correction{Boundary conditions and time integration are not addressed.}{Reviewer 1, Remark 3: semidiscrete, shocks} 
  The novelty is that it combines arbitrary order of convergence,  orthogonal and non-orthogonal structured curvilinear staggered meshes, and the applicability to a wide variety of continuous operators, involving chain rules and nonlinear advection, as illustrated by the shallow-water equations. 
  Experiments show exact conservation and convergence corresponding to expected order.
\end{abstract}

\textbf{Symmetry-preserving discretizations, Mimetic methods, Finite-difference methods, Mass, momentum and energy conservation, Curvilinear staggered grid}

\section{Introduction and motivation}
   Computer simulations of a flow phenomenon require the \discretization\ of the flow properties, reducing the number of values needed to represent the flow state from infinite to some large finite number.  In the resulting discrete model of the flow phenomenon, the (continuous) differential operators have been replaced by (discrete) difference operators.  Unfortunately, not all properties of the differential operators are automatically inherited by their discrete approximations.  The chain and product rules needed in the manipulation of nonlinear equations, for example, do not always work in discrete cases.  Moreover, symmetry and positiveness may be lost in the \discretization\ process, mass, momentum, and energy may not be conserved, aliasing errors can occur, and  duality and self-adjointness of the differential operators may be violated \cite{Bla96SQ, Kok09, Lip14MS}.

   Symmetry-preserving methods, or mimetic methods, aim to preserve certain properties of a continuous operator in its discrete counterpart \cite{Kok09}.  The use of symmetry-preserving \discretizations\ makes it possible to construct discrete models which allow all the manipulations needed to prove stability and (discrete) conservation in the same way they were proven in the original continuous model. In this paper, we will present a new finite-difference symmetry-preserving discretization of arbitrary order on structured curvilinear staggered grids.

   There are a variety of symmetry-preserving \discretizations\ available in the literature. In \cite{Hof12V}, an exhaustive overview is given of different techniques to obtain mass- or energy-conserving methods.  Typically, symmetry properties of differential operators are only automatically preserved in central-difference approximations on uniform, rectilinear grids \cite{Kok09}. Finite-volume methods can be used to construct conservative discretizations for mass and momentum, but it is in general not possible to also obtain energy conservation \cite{Hof12V}.

   In \cite{Ver98V, Ver03V}, a fourth-order symmetry-preserving finite-volume method is constructed using Rich\-ardson extrapolation of a second-order symmetry-preserving method \cite{Vel92R}. The extension to  unstructured collocated meshes is presented in \cite{Tri14LOPSV}, and an application can be seen in \cite{Ver06V}.  The extension to upwind discretizations was made in \cite{Vel08L}, and a discretization for the advection operator for curvilinear collocated meshes was found in \cite{Kok06}. In \cite{Kok09}, the method is extended to non-uniform structured curvilinear collocated grids by deriving a discrete product rule. Furthermore, a symmetry-preserving method that conserves mass and energy for compressible-flow equations with a state equation is described in \cite{Hof12V}. For rectilinear grids, this method works well, but it is challenging to let this method work for unstructured grids \correction{\cite{Hof12V}}{Reviewer 1, Remark 5: Literature}. Finally, in \cite{Roz14VKV}, a symmetry-preserving discretization for curvilinear collocated and rectilinear staggered meshes is found exploiting the skew-symmetric nature of the advection operator on square-root variables.

   Another option to preserve symmetry is to use discrete filters to regularize the convective terms of the equation \cite{Tri11V, Leh12BRPSO}. The combination of a symmetry-preserving discretization and regularization for compressible flows is studied in \cite{Roz15VKV}.

   Mimetic finite-difference methods also mimic the important properties of differential operators. An interesting review is given in \cite{Lip14MS}, and recently, a second-order mimetic discretization of the Navier-Stokes equations conserving mass, momentum, and kinetic energy was presented in \cite{Oud16HVH}. 
   
   Castillo et al. have provided a framework for mimetic operators in \cite{Cas03G, Cas13M}. A second- and fourth-order mimetic approach is constructed for non-uniform rectilinear staggered meshes in \cite{Bat09C, Bla16RCGC, Cor18C}. In \cite{Pue14FHCC} their method is extended to curvilinear staggered meshes, but discrete conservation of mass, momentum and energy is not shown. A second-order mimetic finite-difference method for rectilinear staggered meshes is also constructed in \cite{Sol17GGRO}.
   
   Other mimetic finite-difference methods use algebraic topology to design and analyze compatible discrete operators corresponding to a continuous formulation \cite{Boc06H, Kre11PG, Rob11S}. 
In order to construct a discrete de Rham complex, certain conditions on reconstruction and reduction operators are imposed: they should be conforming, which means that the reconstruction is a right inverse of the reduction \cite{Boc06H}, they should be constant preserving \cite{Bre10B}, and the interpolation operator should commute with the differential operator \cite{Bre10B}. In \cite{Bre14BM} a nice overview of mimetic methods is given. Discrete exterior calculus (DEC) is also related to these mimetic approaches \cite{Hir03}. However, in the papers about these mimetic methods, the discrete conservation of mass, momentum and energy is not studied. In \cite{Pal17G}, a mass, energy, enstrophy and vorticity conserving method is given for unstructured finite-element meshes, using also a conserving time integrator. Note that the method is not applied to staggered meshes, conservation of momentum is not mentioned, and that only incompressible models are used (chain rules are not needed and the advection operator is easier to process). 
More about mimetic time integration is found in \cite{Ste16}.
   
   Another class of symmetry-preserving methods uses the DG method \cite{Win17WGK, Gas16WK, Del14HZ}. In these papers, mass, momentum and \correction{total energy}{Reviewer 2, Remark 7: energy equation}  are conserved on curvilinear meshes, but staggering is not applied. The method is extended to shock capturing and positivity preservation in \cite{Win18WGW}.

   All the existing models in the literature have their own advantages and disadvantages. In general, they are not at the same time applicable for arbitrary discretization orders, sophisticated operators such as the advection operator, or the chain rule, and curvilinear staggered grids. The current paper presents a new discretization method that can handle these requirements simultaneously.

   In this work, we first introduce some concepts by using a Galerkin-type approach, which is closely related to existing mimetic methods, and was also studied in \cite{Hof17V}. Then, we present our new symmetry-preserving finite-difference technique for discretization in space. The novelty of this work is that it combines several important requirements for discretizations: the symmetry-preserving discretization is made for arbitrary order of accuracy; the method works for orthogonal and non-orthogonal structured curvilinear staggered meshes; and the method can be applied to a wide variety of continuous operators, involving chain rules and nonlinear advection, as will be illustrated by the shallow-water equations. The experiments show exact conservation of mass, momentum and energy, and convergence of the approximations corresponding to the expected order. 
   \correction{The approach is very similar to the Richardson extrapolation scheme used in \cite{Kok06} and \cite{Kok09}, but the current approach leads to smaller stencils, especially for high orders of accuracy.}{Reviewer 1, Remark 5: Literature} 
   \correction{Apart from uniform, orthogonal grids, the experiments also use nonorthogonal curvilinear grids, in which the angles between grid lines are as small as 15${}^{\mbox{\scriptsize o}}$.}{Reviewer 2, Remark 11: 3D and non-orthogonal grids}
       \correction{
       The energy equation, derived from the continuity, momentum and state equations, does not have its own discretization. 
       Instead, discrete energy conservation is derived using the symmetries of the discrete operators \cite{Hof19V2}.
       One of the symmetries of interest in this paper is that the adjoint of the gradient is minus the divergence. See Table \ref{Tab: properties}, for more details.} {Reviewer 2, Remark 4: announce properties early}

   \correction{The subject of this paper is symmetry-preserving space discretization. In order to focus on this topic, issues concerning boundary conditions and time integration are not addressed. Instead, a periodic domain is used, and a standard time-integration method is used with such a small time step, that the time-integration errors are negligible.}{Reviewer 1, Remark 3: semidiscrete, shocks} 
   \correction{
   Since the shallow-water equations are the main area of interest for the authors, all examples presented in this paper will be 2D, although the method can also be applied in 1D or 3D.}{Reviewer 2, Remarks 1, 11 and 12:3D}
   
   The outline of this paper is as follows: in Section \ref{Sec: 4 steps}, we present the three different models that will be used in the rest of the paper, and in Section \ref{sec:curvilinear}, we give some information about the curvilinear grids we use. Section \ref{Sec: totals} contains the desired discretization properties, and the new symmetry-preserving discretization is explained in Section \ref{sec:numerical}. The effectivity of this new method is investigated for the different models in Section \ref{sec:results}. We conclude with a discussion of our method and future work in Section \ref{sec:conclusion}.

\section{Models}
\label{Sec: 4 steps}
   \correction{
   In this paper, three different models are used to explain and test the symmetry-preserving ideas. 
   Each of these models is presented in this section with respect to the following aspects:
   \begin{itemize}
   \item  The model equations in terms of continuity, momentum, and state equations;
   \item  A compact representation of the discrete model, including discrete operators like ${\sf DIV}$, the discrete divergence, and ${\sf GRAD}$, the discrete gradient;
   \item  \correction{A consistent energy equation, derived from the continuity, momentum and state equations \cite{Hof19V2}.}{Reviewer 2, Remark 7: energy equation}

    The \correction{energy density}{Reviewer 2, Remark 7: energy equation} is given by the sum of the kinetic and internal energy \correction{densities}{Reviewer 2, Remarks 7: energy equation}: $e = e_{kin}+e_{int}$ (see Table \ref{Tab: mass, momentum and energy} for the definition of energy \correction{density}{Reviewer 2, Remark 7: energy equation} in each model). 
   The energy equation will be used to show energy conservation of the models.
   \end{itemize}

   Section \ref{Sec: totals} discusses the properties that the discrete operators ${\sf DIV}$, ${\sf GRAD}$ and others are expected to have, for the discrete models to conserve mass, momentum and energy, after which
   Section \ref{sec:numerical} discusses how such operators can be constructed.}{Reviewer 2, Remark 5: announce goals in S 2}

   \begin{itemize}

      \item {\bf Linear-wave equations}

         The simplest equation that can be used to discuss symmetry preservation, is the linear-wave equation, in which the evolution of the pressure $p$ \correction{in a domain $V$,}{Reviewer 1, Remark 1 and Reviewer 2, Remark 6: Introducing quantities} 
         the flow velocity $\vec v$ and the density $\rho$, are given by the continuity,  momentum, and state equations:
         \begin{eqnarray}
            \paf \rho t+ \rho_0 \nabla \cdot \vec v = 0, & \displaystyle
            \paf{\vec v} t + \frac 1{\rho_0}\nabla p = 0, & p = c^2 \rho,
         \label{eqn:linear-wave equations}
         \end{eqnarray}
         \correction{where $\rho_0$ is a constant 'reference' density, and $c$ is the speed of sound (i.e. the propagation speed of waves).}{Reviewer 1, Remark 1 and Reviewer 2, Remark 6: Introducing quantities}  Initial conditions are specified at \correction{initial time $t=0$}{Reviewer 1, Remark 1 and Reviewer 2, Remark 6: Introducing quantities}.

         The discrete linear-wave equations are given by
         \begin{eqnarray}
            \naf{\,\tt rho} t + \rho_0 ~{\sf DIV}~{\tt v} = 0, & \displaystyle
            \naf{\,\tt v} t   + \frac 1{\rho_0}{\sf GRAD}~{\tt p} = 0, & {\tt p} = c^2~{\tt rho},
            \nonumber
         \end{eqnarray}
         where {\tt p}, {\tt rho} and {\tt v} are the vectors with the discrete pressures, densities and velocities, and where {\sf DIV} is the discrete divergence and {\sf GRAD} is the discrete gradient.

         The continuity, momentum and state equations can be combined into the following energy equation:
         \[ \paf et + \nabla \cdot p \vec v = 0,  \]
         where $e$ is the energy \correction{density}{Reviewer 2, Remark 7: energy equation}.

         This model gives us the opportunity to introduce the relation between symmetry preservation and conservation, as well as the curvilinear staggered grid and staggered velocity components. Two approaches are used for the construction of a discretization. First, we use a Galerkin-type approach similar to the one used in \cite{Hof17V}, and secondly, we apply a finite-difference approach.

      \item {\bf Compressible-wave equations}

         We introduce a non-linearity into the system by including density variations in the continuity equation, and so obtain com\-press\-i\-ble-wave equations (without an advection term):
         \begin{eqnarray}
            \paf \rho t+ \nabla \cdot \rho \vec v = 0, & \displaystyle
            \paf{\vec v} t + \nabla Q(p) = \paf{\vec v} t + \frac 1{\rho}\nabla p = 0, & \rho = R(p),
            \nonumber
         \end{eqnarray}
         \correction{where the pressure-dependent density function $R$ is assumed positive, continuous and monotonically non-decreasing, and where $Q$ is given by $Q(p):=\int^p 1/R(y)~\mbox dy$}{Reviewer 1, Remark 1 and Reviewer 2, Remark 6: Introducing quantities}. The discrete com\-press\-i\-ble-wave equations are given by
         \begin{equation}
            \naf{\,\tt rho} t + {\sf DIV}\tilde {\sf r}~ {\tt v} = 0, \quad
                 \naf{\,\tt v} t + {\sf GRAD}~Q({\tt p}) = 0, \quad {\tt rho} = R({\tt p}),
         \label{eq: compressible-wave equations, discrete}
         \end{equation}
         where the operator ${\sf DIV}\tilde{\sf r}$ is a discretization of the operator $(\nabla \cdot \rho)$: this means that $({\sf DIV}\tilde{\sf r} \ {\tt v})_i$ approximates $(\nabla \cdot \rho \vec v)(\vec x_i)$.  In Section \ref{sec:rGRAD}, the explicit construction of the operator ${\sf DIV}\tilde{\sf r}$ will be discussed.

         The corresponding energy equation is
         \[ \paf et + \nabla \cdot \left( e_{int} + \rho_0 \int^p \frac {1}{R(q)} \ \mbox dq \right) \vec v = 0.  \]

         The model is analyzed for an arbitrary state equation $\rho=R(p)$, but tests are only conducted for the state equation $p = c^2\rho$, because then an exact solution is available to verify the results.

      \item {\bf Isentropic compressible Euler equations}

         Finally, after introducing an advection term, a symmetry-preserving discretization is formed for the following equations of isentropic compressible Euler gas dynamics \cite{Lev02}:
         \begin{eqnarray}
            \paf \rho t+ \nabla \cdot \rho \vec v = 0, & \displaystyle
            \paf{} t \rho \vec v + \nabla\cdot(\rho \vec v\otimes \vec v) + \nabla p = 0, & \rho = R(p).
            \nonumber
         \end{eqnarray}
         The discrete system is given by
         \begin{equation}
            \naf{\,\tt rho} t + {\sf DIVr}~ {\tt v} = 0, \
            \naf{\,\tt rv} t + {\sf ADVEC}~{\tt v} + {\sf GRAD}~{\tt p} = 0, \ {\tt rho} = R({\tt p}).
         \label{eqn:Euler}
         \end{equation}
         The vector {\tt rv} contains local momentum values, and is defined by ${\tt rv} := \diag({\sf Interp}_{v\leftarrow c} {\tt rho}) {\tt v}$. ${\sf Interp}_{v\leftarrow c}$ indicates interpolation from the pressure grid points to the velocity grid points in the staggered grid. The operator ${\sf DIVr}$ is similar to the operator ${\sf DIV}\tilde{\sf r}$: it also approximates the operator $(\nabla \cdot \rho)$, but is constructed in a different way (Section \ref{sec:Advection}). \correction{Section \ref{sec:Advection} also shows the construction of the discrete advection operator {\sf ADVEC}.}{Reviewer 1, Remark 1 and Reviewer 2, Remark 6: Introducing quantities}

         The following energy equation can be derived from the continuity, momentum and state equations of the isentropic compressible Euler equations:
         \[ \paf et + \nabla \cdot (e+p) \vec v = 0. \]

         \correction{Choosing the state equation $\rho = R(p) = \sqrt{2p/g}$ changes these equations into the shallow-water equations \cite{Lev02}, which are investigated in this paper.}{Reviewer 1, Remark 1 and Reviewer 2, Remark 6: Introducing quantities} 
   \end{itemize}

   For each model, an exact solution is created using the propagation speeds and Riemann invariants as given in Table \ref{Tab: speeds and invariants}.
   \correction{For the compressible-wave equations and the isentropic compressible Euler
   equations, the exact solutions develop shocks after some time.  After that
   moment, the (weak) solution becomes discontinuous, and although mass and
   momentum are still conserved, energy is not.  
   In Section \ref{sec:results}, we 
   \ignore{Tables \ref{tab:conservation_compressible} and
   \ref{tab:conservation_shallow}} show that the symmetry-preserving method
   conserves mass, momentum and energy until the moment the shock occurs, 
   and that the method converges with the expected order at a time before the shock occurs. 
   At the moment the shock occurs, the convergence is slower.
   After the appearance of the shock, the discrete solution, which still conserves energy, cannot possibly be an accurate approximation of the correct (weak) solution, in which energy is lost.  The comparison of the solutions therefore stops when the shock appears.
   The method could be extended with artificial
   viscosity (for instance in the form of flux limiters), which dissipate some
   of the energy in the system, and allow correct approximation of the
   discontinuous solution, but these techniques are beyond the scope of this article.

   An example with a developing shock was chosen to show the properties of the method, because it shows that the method is accurate when it can be expected to, that it always conserves mass, momentum and energy, and that there are cases (discontinuous solutions) where the method needs an extension before it can be used to calculate meaningful solutions. All these aspects of the method would not be illustrated by an example where shocks do not develop.}{Reviewer 1, Remark 3: semidiscrete, shocks} 
   \begin{table}[ht!]
      \centering
      \caption{Propagation speeds and Riemann invariants used to construct exact one-dimensional solutions for the three models \cite{Lev02}. The one-dimensional velocity is given by $v$.}
      \resizebox{\textwidth}{!}{% use resizebox with textwidth
         \begin{tabular}{l|l|l|l|l}
                                          & Linear-wave eq. & Compressible-wave eq. & Shallow-water eq. \\
            \hline
                                          &                 &                       &                   \\
            {\bf Propagation speed $V$}   &                 &                       &                   \\
            ~~~{forward} $V_+$            &     $c$         &$\displaystyle \frac 12
                                                             (v+\sqrt{v^2+4c^2})$   & $v+\sqrt{g\rho}$  \\[2ex]
                                          &                 &                       &                   \\
            ~~~{backward} $V_-$           &    $-c$         &$\displaystyle \frac 12
                                                               (v-\sqrt{v^2+4c^2})$ & $v-\sqrt{g\rho}$  \\[2ex]
               \hline
                                          &                 &                       &                   \\
            {\bf Riemann invariant $F$}   &                 &                       &                   \\
            ~~~{forward} $F_+$            
                                          &$\rho_0 v+c\rho$ & $\displaystyle
                                                              \ln(\rho V_+)-
                                                              \frac{vV_-}{2c^2}$    & $v+2\sqrt{g\rho}$ \\[2ex]
                                          &                 &                       &                   \\
            ~~~{backward} $F_-$           
                                          &$\rho_0 v-c\rho$ & $\displaystyle \ln\left(
                                                                                       \frac{\rho}{V_+}\right)
                                                                                       -\frac{vV_+}{2c^2}$ & $v-2\sqrt{g\rho}$ \\[2ex]
         \end{tabular}
      }
   \label{Tab: speeds and invariants}
   \end{table}

\section{Curvilinear grid}\label{sec:curvilinear}
   In this section, we introduce a parametrization for curvilinear grids. This is first done for a collocated mesh (used for scalar fields) and then for a staggered grid (used in vector fields).
   
   \subsection{Collocated grids and scalar fields}
   We introduce a uniform grid in computational-grid space, $\vec \xi_{c,i}$, that satisfies
   \[ \vec \xi_{c,i+m_xj+m_xm_yk} = (i\Delta\xi, \ j\Delta\eta, \ k\Delta\zeta)^\top, \]
   where $(m_x,m_y,m_z)$ relate to the number of cells in each direction, and $(\Delta \xi, \Delta \eta, \Delta \zeta)$ are the corresponding cell widths. The subscript $c$ is chosen because the grid points can be seen as the cell centers of control volumes with vertices  $\vec \xi_{c,i} + (\pm\Delta\xi, \ \pm\Delta\eta, \ \pm\Delta\zeta)^\top/2$.

   The map ${\vec X}$ relates this uniform grid $\{\vec \xi_{c,0},\ldots,\vec \xi_{c,m_xm_ym_z-1}\}$ in computational-grid space to a curvilinear grid $\{\vec x_{c,0},\ldots,\vec x_{c,m_xm_ym_z-1}\}$ in physical space:
   \[ \vec x_{c,i} = \vec X(\vec \xi_{c,i}).\]
   The grid formed by the points $\xi_{c,i}$ is called the {\em pressure grid}, because this is the grid used to sample scalar fields, such as the pressure.

   \correction{The properties of the mapping $\vec X$, along with those of the solution, determine the convergence rate of the discretization used. For clarity of presentation, we have only used smooth ($C^\infty$) mappings.}{Reviewer 1, Remark 1 and Reviewer 2, Remark 6: Introducing quantities} 

   \subsection{Staggered grids and vector fields}
\begin{figure}[ht!]
       \centering
        \begin{tikzpicture}
       \draw (0,-0.5)--(0,2.5);
       \draw (3,2.5)--(3,5.5);
       \draw (-0.5,-0.5)--(3.5,3.5);
       \draw (-0.5,1.5) -- (3.5,5.5);
       \draw[fill=black] (1.5,2.5) circle (0.5mm);
       \node at (1.5,2.2) {$\vec x_{c,i}$};
       \draw[fill=black] (0,1) circle (0.5mm);
       \node at (-0.5,1) {$\vec x_{e,i-1}$};
       \draw[fill=black] (3,4) circle (0.5mm);
       \node at (3.5,4) {$\vec x_{e,i}$};
       \draw[fill=black] (1.5,1.5) circle (0.5mm);
       \node at (1.8,1) {$\vec x_{n,i-m_x}$};
       \draw[fill=black] (1.5,3.5) circle (0.5mm);
       \node at (1.5,4) {$\vec x_{n,i}$};
       \end{tikzpicture}
       \caption{Numbering of the grid points for a staggered grid in 2D. The points $\vec x_{c,i}$ belong to the pressure grid cells (cell centers), and $\vec x_{e,i}$, $\vec x_{n,i}$ belong to the velocity grid, with grid points located at the east and north cell faces, respectively.}\label{fig:numbering}
      \end{figure}
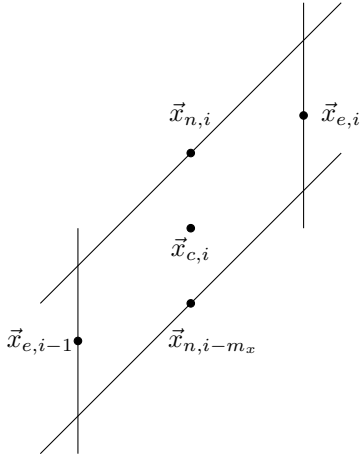

      In a staggered grid, not only the pressure grid is used, but also a velocity grid. The velocity grid consists of grid points that are shifted by half a grid space, and therefore correspond to {\em cell-face centers}:
       \begin{eqnarray}
          \vec x_{e,i} := \vec X\left(\vec \xi_{e,i}\right) :=
                          \vec X\left(\vec \xi_{c,i} + \frac 12 {\mbox{\scriptsize $\left(\!\!\begin{array}l \Delta \xi\\0\\0\end{array}\!\!\right)$}}\right),&&
          \vec x_{n,i} := \vec X\left(\vec \xi_{n,i}\right) :=
                          \vec X\left(\vec \xi_{c,i} + \frac 12 {\mbox{\scriptsize $\left(\!\!\begin{array}l 0 \\ \Delta \eta\\0\end{array}\!\!\right)$}}\right),
\nonumber \\
          \vec x_{t,i} := \vec X\left(\vec \xi_{t,i}\right) :=
                          \vec X\left(\vec \xi_{c,i} + \frac 12 {\mbox{\scriptsize $\left(\!\!\begin{array}l 0 \\ 0 \\ \Delta \zeta\end{array}\!\!\right)$}}\right),
      \nonumber
      \end{eqnarray}
      where the subscripts $e$, $n$ and $t$ stand for 'east', 'north' and 'top', respectively.
      The numbering of the grid points is visualized in Figure \ref{fig:numbering}.
      
      To find vector-field discretizations in curvilinear spaces, we should first introduce the notation of {\em grid-aligned velocity components}. A vector field $\vec v$ is represented by three scalar functions $v_x$, $v_y$ and $v_z$, using a {\em local grid orientation} $(\vec r_x, \vec r_y, \vec r_z)$,
      \begin{equation}
         \vec v(\vec x) = v_x(\vec x) \vec r_x(\vec x)  + v_y(\vec x) \vec r_y(\vec x)  + v_z(\vec x) \vec r_z(\vec x),
      \label{eq: v from components}
      \end{equation}
      where $\vec r_x$, $\vec r_y$ and $\vec r_z$ are orthonormal. The scalar functions $v_x$, $v_y$ and $v_z$ can be calculated from the vector field $\vec v$ using inner products:
      \begin{eqnarray}
      v_x(\vec x) = \vec r_x(\vec x)
         \cdot \vec v(\vec x),&
      v_y(\vec x) = \vec r_y(\vec x)
         \cdot \vec v(\vec x),&
      v_z(\vec x) = \vec r_z(\vec x)
         \cdot \vec v(\vec x).
      \end{eqnarray}

\begin{figure}[ht!]
         \centering
         \includegraphics[scale=0.5]{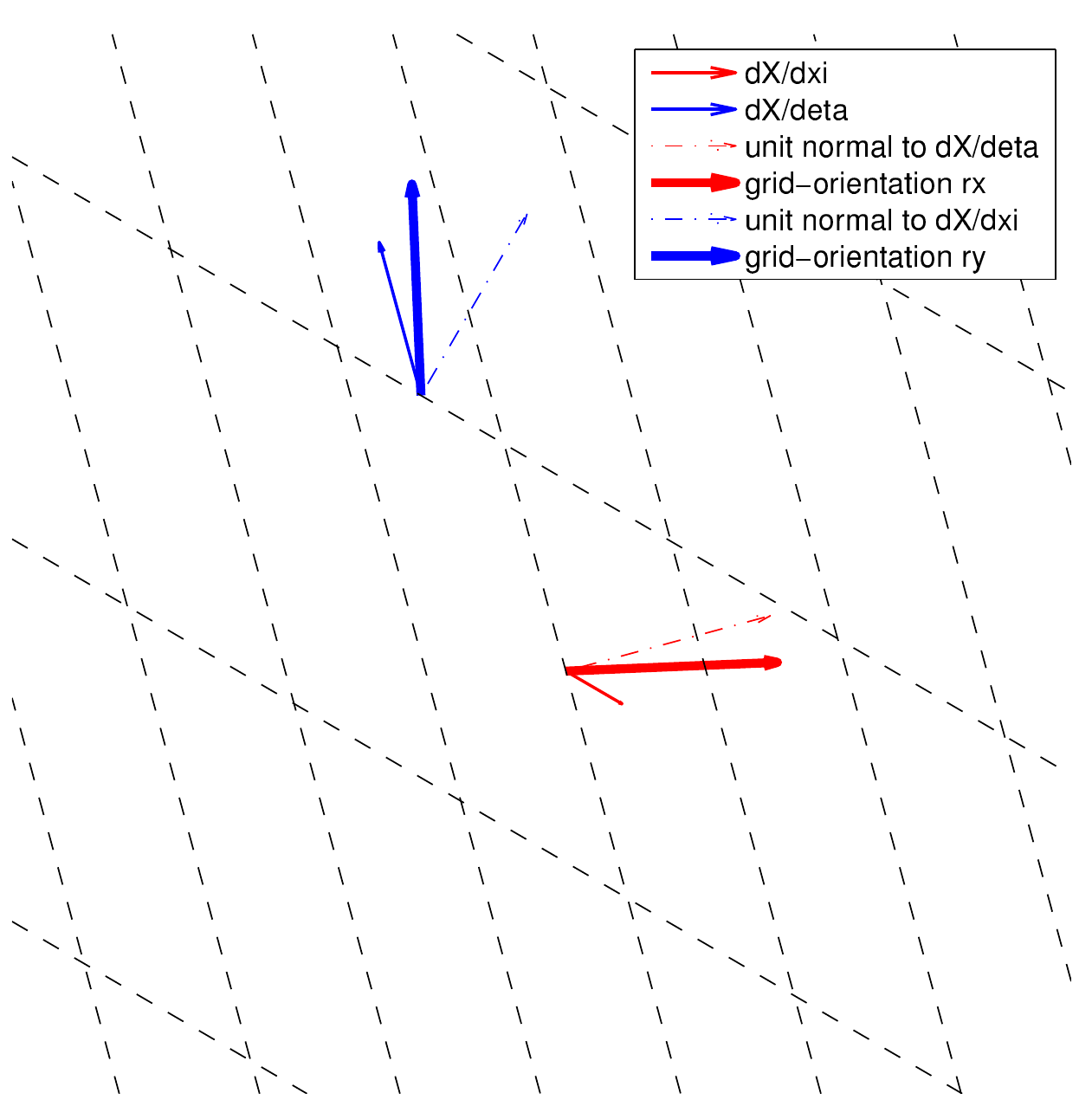}
         \caption{Example of a 2D curvilinear grid and corresponding relations between the Cartesian grid $(x,y)$ and $(\xi, \eta)$. The solid thin lines correspond to the orientation of the grid lines. The dashed lines are orthogonal to these grid lines. The solid thick lines form the local grid orientation, computed using the SVD technique. Note that $\vec r_x$ is closer to the normal (dashed red line) than to the grid line (thin red line), and $\vec r_y$ is farther from the normal. This is because the cells are stretched in the $y$-direction. The vectors are located at the staggered velocity grid points $\vec{x}_e$ and $\vec{x}_n$.}
      \label{fig:plaatje}
      \end{figure}

      \correction{Typically \cite{Bee95NW, Wes01}, a combination of {\em covariant}
      directions (i.e. parallel to grid lines) and {\em contravariant}
      directions (i.e. perpendicular to two of the grid lines) is used to
      construct a local grid orientation. Often \cite{Bee95NW}, this
      leads to a different grid orientation in $\xi$-cell faces than in $\eta$-
      or $\zeta$-cell faces. 
      For the method in this paper, a different kind of grid orientation is needed, because a simple representation of the kinetic energy requires a local grid orientation that is
      \begin{itemize}
         \item available in all points in space, not just in grid points;
         \item is the same for all quantities and equations;
         \item consists of three orthogonal directions.
      \end{itemize}

      In uniform Cartesian grids, the covariant and contravariant directions are the same, and the local grid orientation is given by $\vec r_x=(1,0,0)$, $\vec r_y=(0,1,0)$ and $\vec r_z=(0,0,1)$.  In the general case, the contravariant and covariant directions are not the same, and an orthogonal grid orientation cannot be made from them.

      An orthogonal, grid-aligned basis for vectors is based on the singular-value decomposition (SVD) of the grid directions, given by the vector ${\bf h}$ with the singular values and unitary matrices ${\bf P}$ and ${\bf Q}$, such that
      \begin{eqnarray}
         \left( \begin{array}{lll}
            \displaystyle \paf{\vec X}\xi  &
            \displaystyle \paf{\vec X}\eta &
            \displaystyle \paf{\vec X}\zeta
         \end{array}\right)
         = {\bf P} ~\diag({\bf h}) ~{\bf Q}.
      \end{eqnarray}
      The orientation ($\vec r_x, \vec r_y, \vec r_z$) is now found by using only the unitary rotation matrices of the singular-value decomposition:
      \[ [\vec r_x,\vec r_y, \vec r_z] := {\bf P} ~{\bf Q}. \]
      This local-orientation matrix is exactly orthogonal, and forms a compromise between the covariant and contravariant directions. A 2D example of the computation of such local orientation is given in Figure \ref{fig:plaatje}.}{Reviewer 1, Remark 2, SVD}

      A discrete vector field is represented by a vector with values for each direction: ${\tt v}=({\tt vx}^\top \hspace{-0.2cm},{\tt vy}^\top\hspace{-0.2cm},{\tt vz}^\top)^\top$ contains components ${\tt vx}$, {\tt vy}, and ${\tt vz}$ that are located at the grid points $\vec{x}_e$, $\vec{x}_n$, and $\vec{x}_t$, respectively. The components are given by
      \begin{equation}
          {\tt vx}_i = \vec r_x(\vec x_{e,i}) \cdot \vec v(\vec x_{e,i}), \quad
          {\tt vy}_i = \vec r_y(\vec x_{n,i}) \cdot \vec v(\vec x_{n,i}), \quad
          {\tt vz}_i = \vec r_z(\vec x_{t,i}) \cdot \vec v(\vec x_{t,i}).
          \label{eqn:discretevectorfield}
      \end{equation}

      Discrete samplings of the local grid orientation are stored in vectors ${\tt r**\_at\_*}$. As an example, we give the expressions for $\vec r_x$ at the cell centers or east cell-face centers, and for $\vec r_y$ at the north cell-face centers:
      \begin{eqnarray}
          \vec r_x(\vec x_{c,i}) = \left(\begin{array}{c}
                                 {\tt rxx\_at\_c}_i \\
                                 {\tt rxy\_at\_c}_i \\
                                 {\tt rxz\_at\_c}_i \end{array}\right),&
          \vec r_x(\vec x_{e,i}) = \left(\begin{array}{c}
                                 {\tt rxx\_at\_e}_i \\
                                 {\tt rxy\_at\_e}_i \\
                                 {\tt rxz\_at\_e}_i \end{array}\right),&
          \vec r_y(\vec x_{n,i}) = \left(\begin{array}{c}
                                 {\tt ryx\_at\_n}_i \\
                                 {\tt ryy\_at\_n}_i \\
                                 {\tt ryz\_at\_n}_i \end{array}\right).
      \end{eqnarray}

      Unlike in a uniform Cartesian grid, the discrete representation of the constant vector field $(1,0,0)$ is not a constant vector, and neither does it have zeros in the $n$-points and $t$-points of the grid.  To compute the discrete representations {\tt c100}, {\tt c010} and {\tt c001} of the constant fields 
      \[
       \vec v(\vec x) = \vec c_{(1,0,0)}(\vec x) = (1,0,0), \quad
       \vec v(\vec x) = \vec c_{(0,1,0)}(\vec x) = (0,1,0), \quad
       \vec v(\vec x) = \vec c_{(0,0,1)}(\vec x) = (0,0,1),
      \]
      we use equation (\ref{eqn:discretevectorfield}) to find:
      \begin{equation}
      {\tt c100} = 
      \begin{pmatrix}
       {\tt rxx\_at\_e} \\
       {\tt ryx\_at\_n} \\
       {\tt rzx\_at\_t}
      \end{pmatrix}, \quad
      {\tt c010} =
      \begin{pmatrix}
       {\tt rxy\_at\_e} \\
       {\tt ryy\_at\_n} \\
       {\tt rzy\_at\_t}
      \end{pmatrix}, \quad
      {\tt c001} = 
      \begin{pmatrix}
       {\tt rxz\_at\_e} \\
       {\tt ryz\_at\_n} \\
       {\tt rzz\_at\_t}
      \end{pmatrix}. \label{eqn:c100}
      \end{equation}

\section{Desired discretization properties} \label{Sec: totals}
   In this section, we investigate the properties that a symmetry-preserving discretization should satisfy. Therefore, we first introduce the different scalar products used. 
   
   \subsection{Scalar products}
   Quantities such as densities and pressures form scalar fields, which are discretely approximated by vectors with a value for each pressure grid point. Let $a$ and $b$ be two continuous scalar fields, that correspond to the discrete scalar fields ${\tt a}$ and ${\tt b}$. The integral $\int_V ab \ \mbox dV$ can be approximated by the following scalar product:
   \begin{equation}
      \langle {\tt a}, {\tt b} \rangle_c := {\tt a}^\top \diag({\tt dVc}) {\tt b}, \label{eqn:definitionvectorproduct}
   \end{equation}
   where only real numbers are used, and the vector ${\tt dVc}$ holds the integration weights for each pressure grid point. For instance, the discrete mass, {\tt M}, is calculated from the vector ${\tt rho}$ of discrete densities according to ${\tt M} :=  \langle{\tt c1} , {\tt rho}\rangle_c$, where ${\tt c1}$ is the (constant) vector of only ones.

   \bigskip
   The situation is more complicated when vector fields are involved. Since the grid orientations have been chosen orthogonal, the scalar product $\langle{\tt v},{\tt w}\rangle_v$ of two discrete vector fields ${\tt v}$ and ${\tt w}$ can be written as the sum of the scalar products of the components:
   \begin{align}
      \langle{\tt v},{\tt w}\rangle_v &:=
      \langle{\tt vx} , {\tt wx}\rangle_e + \langle{\tt vy} , {\tt wy}\rangle_n + \langle{\tt vz} , {\tt wz}\rangle_t
      \nonumber \\
       &:= 
       {\tt vx}^\top ~\diag({\tt dVe})~{\tt wx} + {\tt vy}^\top ~\diag({\tt dVn})~{\tt wy}+ {\tt vz}^\top ~\diag({\tt dVt})~{\tt wz},
   \end{align}
   where ${\tt dVe}$, ${\tt dVn}$ and ${\tt dVt}$ are the vectors holding the integration weights for the grid points on the three directions of the velocity grid. 

   As an example, the total momentum $\vec{\tt M}$ is calculated for the compressible-wave model, using a discrete vector field {\tt rv} of local momentum values: 
   \[
     \vec {\tt M} :=
      \left(
         \langle{\tt c100} , {\tt rv}\rangle_v, \ \langle{\tt c010} , {\tt rv}\rangle_v, \ \langle{\tt c001} , {\tt rv}\rangle_v
      \right)^\top.
   \]
     
   \bigskip
   The scalar products introduced in this section determine the corresponding adjoints. For example, for a matrix ${\sf A}$ that maps values on pressure grid points to values on $e$-points, the adjoint ${\sf A}^*$ is given by
   \[
    {\sf A}^* = \diag({\tt dVc})^{-1}~{\sf A}^\top~\diag({\tt dVe}),
   \]
   since in that case, we have for all vectors ${\tt x}$ on $e$-points and ${\tt y}$ on pressure points:

   \[     \langle{\tt x}, {\sf A}~{\tt y}\rangle_e = {\tt x}^\top \diag({\tt dVe}) {\sf A}~{\tt y} = {\tt x}^\top ({\sf A}^*)^\top \diag({\tt dVc})~{\tt y} = \langle{\sf A}^*{\tt x}, {\tt y}\rangle_c. \]
   The scalar products and adjoints of this section are used to define the properties that discrete operators should satisfy.

   \subsection{Operator properties}\label{sec:properties}
   In this section, we discuss Table \ref{Tab: properties}, that shows some properties of the operators used to define the models from Section \ref{Sec: 4 steps}, along with their discrete equivalents. In \cite{Roz14VKV}, symmetry properties are used instead of the null space properties that are given in the table. These symmetry properties are applied to curvilinear collocated grids, and therefore this approach could not be applied to curvilinear staggered grids without the modifications presented in this paper.

   The continuous versions of the properties often contain a boundary integral. For example, the change in total mass in the compressible-wave equations equals
   \begin{equation}
     \paf Mt = \paf{}{t} \int_V \rho \ \mbox dV = - \int_V \nabla \cdot \rho \vec v \ \mbox dV = - \oint_{\delta V} \rho \vec v \cdot \vec n \ \mbox dS. \label{eqn:int}
   \end{equation}
   In this paper, we do not focus on boundary conditions, and therefore, we take periodic domains, such that equation (\ref{eqn:int}) equals zero. We could also use closed walls, which satisfy $\vec v \cdot \vec n = 0$. For other boundary conditions, the total mass, momentum and/or energy may change over time, due to an in- or outflux across the domain boundaries. This is beyond the scope of this paper.

   Many discrete equalities in the table are of the same form, containing an adjoint operator and a vector consisting of ones. As an example, the equality for ${\sf DIV}^* {\tt c1}$ should be read as follows:
         \begin{align*}
            {\sf DIV}^* {\tt c1}=0 &\iff {\tt f}^\top \ \diag({\tt dVc}) \ {\sf DIV}^* {\tt c1}=0 \quad \forall {\tt f}  \\
            &\iff {\tt f}^\top \ \diag({\tt dVc}) \ \left\{ \diag({\tt dVc})^{-1} \ {\sf DIV}^\top \ \diag({\tt dVc}) \right\} \ {\tt c1} = 0 \quad \forall {\tt f} \\
            &\iff {\tt f}^\top \ {\sf DIV}^\top \ \diag({\tt dVc}) \ {\tt c1} = 0 \quad \forall {\tt f} \\
            &\iff {\tt c1}^\top \ \diag({\tt dVc})\ {\sf DIV}  \ {\tt f} = 0 \quad \forall {\tt f} \\
            &\iff \text{The discrete integral of } {\sf DIV} \ {\tt f} \text{ equals zero} \quad \forall {\tt f},
         \end{align*}
   which relates to its continuous equivalent. In Section \ref{sec:numerical}, we will construct the discrete operators that satisfy these properties.
   
      \begin{table}[htp]
      \caption{Several properties of the operators used in the models. The left null space properties are used to prove conservation of mass and momentum. For energy conservation, all properties are required.}
      \resizebox{\textwidth}{!}{% use resizebox with textwidth
      \begin{tabular}{l|l}
      \hline
         \multicolumn{2}{c}                  {\bf Left null space properties}                                          \\
         \multicolumn{1}{c}{\bf Continuous version } & \multicolumn{1}{c}{\bf Discrete equivalent}                   \\
         \hline
                                                     &                                                               \\
         $\displaystyle \int \nabla \cdot \vec f~\mbox dV = \oint_{\delta V} \vec f \cdot \vec n ~\mbox dS = 0$
                                                     & ${\sf DIV}^* {\tt c1}=0$                                      \\[2ex]
         $\displaystyle \int \nabla f ~\mbox dV           = \oint_{\delta V} f \vec n~\mbox dS = \vec 0$
                                                     & ${\sf GRAD}^*~{\tt c100} = {\sf GRAD}^*~{\tt c010}
                                                                                = {\sf GRAD}^*~{\tt c001} = 0$       \\[2ex]
         $\displaystyle \int_V \nabla \cdot( \rho \vec v \otimes \vec f) ~\mbox dV =
                    \oint_{\delta V} \rho \vec f \vec v \cdot \vec n~\mbox dS = \vec 0$
                                                     & ${\sf ADVEC}^*{\tt c1}=0$                                     \\[2ex]
         $\displaystyle \int_V \nabla \cdot \rho \vec f ~\mbox dV =
                    \oint_{\delta V} \rho \vec f \cdot \vec n ~\mbox dS = 0$
                                                     & $ {\sf DIV}\tilde {\sf r}^*~{\tt c1}=0$                       \\[2ex]
      \hline
         \multicolumn{2}{c}          {\bf Zero derivative for constants (right null space properties)} \\
         \multicolumn{1}{c}{\bf Continuous version } & \multicolumn{1}{c}{\bf Discrete equivalent}                   \\
      \hline
                                                     &                                                               \\
         $\displaystyle \nabla \cdot (1,0,0) = \nabla \cdot (0,1,0) =
                           \nabla \cdot (0,0,1) = 0$ &  ${\sf DIV}~{\tt c100} = {\sf DIV}~{\tt c010}
                                                                              = {\sf DIV}~{\tt c001} = 0$            \\[1ex]
         $\displaystyle \nabla 1 = \vec 0$           &  ${\sf GRAD}~{\tt c1}=0$                                      \\[1ex]
      \hline
         \multicolumn{2}{c}             {\bf Chain-rule properties} \\
         \multicolumn{1}{c}{\bf Continuous version } & \multicolumn{1}{c}{\bf Discrete equivalent}                   \\
      \hline
                                                     &                                                               \\
         $\displaystyle R(p) \nabla S(p) =
                                        \nabla Q(p)$ & $\tilde{\sf r}{\sf GRAD}~S({\tt p}) = {\sf GRAD}~Q({\tt p})$  \\[2ex]
         $\displaystyle R(p) \nabla Q(p) = \nabla p$ & ${\sf rGRAD}~Q({\tt p}) = {\sf GRAD}~{\tt p}$                 \\[2ex]
         with $Q(p):=\int^p \frac 1{R(q)}~\mbox dq$ and
         $S(p) :=  \int^p \frac 1{R^2(q)}~\mbox dq$ \\[2ex]
      \hline
         \multicolumn{2}{c}               {\bf Symmetry properties} \\
         \multicolumn{1}{c}{\bf Continuous version } & \multicolumn{1}{c}{\bf Discrete equivalent}                   \\
      \hline
                                                     &                                                               \\
         $\displaystyle \int_V f \nabla \cdot \vec g~\mbox dV  +  \int_V \vec g \cdot \nabla f~\mbox dV =
                   \oint_{\delta V} f\vec g \cdot \vec n ~\mbox dS = 0$
                                                     & ${\sf DIV} + {\sf GRAD}^*=0$                                  \\[2ex]
         $\displaystyle \int_V f \nabla \cdot \rho \vec g~\mbox dV  +  \int_V \vec g \cdot \rho \nabla f~\mbox dV =
                        \oint_{\delta V} \rho f\vec g \cdot \vec n ~\mbox dS = 0$
                                                     & ${\sf DIV}\tilde {\sf r} + \tilde{\sf r}{\sf GRAD}^* = 0$     \\[2ex]
         $\displaystyle \int_V f \nabla \cdot \rho \vec g~\mbox dV  +  \int_V \vec g \cdot \rho \nabla f~\mbox dV =
                        \oint_{\delta V} \rho f\vec g \cdot \vec n ~\mbox dS = 0$
                                                     & ${\sf DIVr} +{\sf rGRAD}^* = 0$                               \\[2ex]
         $\displaystyle \int_V \vec g \cdot \nabla \cdot \rho \vec v \otimes f ~\mbox dV  +
                        \int_V \vec f \cdot \nabla \cdot \rho \vec v \otimes g ~\mbox dV $
                                                     & ${\sf ADVEC} + {\sf ADVEC}^* =
                                                            \diag({\sf Interp}_{v\leftarrow c}~{\sf DIVr}~{\tt v})$  \\[2ex]
         $ = \displaystyle \int_V \vec f \cdot \vec g \nabla \cdot \rho \vec v ~\mbox dV +
                             \oint_{\delta V} \rho (\vec f \cdot \vec g) (\vec v \cdot \vec n)~\mbox dS = \int_V \vec f \cdot \vec g \nabla \cdot \rho \vec v ~\mbox dV$
      \end{tabular}
      }
      \label{Tab: properties}
   \end{table}
   
   \subsection{Conserved quantities}
     In each model, the conserved quantities (total mass, (total momentum) and \correction{total energy}{Reviewer 2, Remark 7: energy equation}) can be identified, which are computed using integrals over the simulation domain $V$. For instance, for the linear-wave equations, the total mass $M:=\int_V\rho~\mbox dV$ is a conserved quantity. Table \ref{Tab: mass, momentum and energy} displays the expressions for each of the conserved quantities in the models used in this paper. The symmetry-preserving discretizations should conserve these quantities in a discrete manner. In this section, we explain how discrete conservation can be proved, thereby using scalar products, and the properties from Table \ref{Tab: properties}. Conservation of discrete mass in the linear-wave equations, for instance, is shown in the following very short proof:
     \[
        \naf{\tt M}t  = \naf{}t \langle{\tt c1}, {\tt rho}\rangle_c =
                           \left\langle{\tt c1}, \naf{~\tt rho}t\right\rangle_c =
                           -\left\langle{\tt c1}, \rho_0 {\sf DIV}~{\tt v}\right\rangle_c =
                           -\left\langle {\sf DIV}^*{\tt c1}, \rho_0~{\tt v}\right\rangle_v = 0.
      \]
      Note the transition from an inner product on the pressure grid to the velocity grid in the second last equality. The velocity $v$ is defined at the cell faces, and ${\sf DIV}\ {\tt v}$ is located at the cell centers. Therefore, the inner product $\left\langle{\tt c1}, \rho_0 {\sf DIV}~{\tt v}\right\rangle_c$ is computed in the pressure grid. Note that $\rho_0$ is constant, and therefore available both on the pressure and on the velocity grid. All proofs for the conservation of mass and momentum follow these same lines, using the left null space properties of the discretizations. 

     The energy-conservation proofs also use the symmetry properties. As an example, we show energy conservation for the isentropic Euler equations.
     \correction{The derivation of the continuous and discrete energy equations for all the models in this paper are presented in \cite{Hof19V2}.}{Reviewer 2, Remark 7: energy equation} 
     In the isentropic Euler equations, energy conservation is shown by calculating the time derivative of the discrete \correction{total energy}{Reviewer 2, Remark 7: energy equation}  {\tt E}, \correction{which is the discrete approximation of the (continuous) \correction{total energy}{Reviewer 2, Remark 7: energy equation}  $E:=\int_V e\mbox dV$}{Reviewer 1, Remark 1 and Reviewer 2, Remark 6: Introducing quantities}:
     \[    \naf{\tt E}t = \naf{}t \left( \langle{\tt c1}, e_{int}({\tt p})\rangle_c + \frac 12 \langle{\tt v},{\tt rv}\rangle_v \right). \]
   If we apply the chain rule twice, we find that the first term equals
   \[ \naf{}t \langle{\tt c1}, e_{int}({\tt p})\rangle_c
          = \left\langle{\tt c1}, \frac{e^\prime_{int}({\tt p})}{R^\prime({\tt p})} \naf{~\tt rho}t \right\rangle_c.  \]
   For the second term, we use the fact that we only consider real numbers, such that for fields ${\tt a}, {\tt b}, {\tt c}$ we have
   \begin{equation}
      \langle {\tt a},\diag({\tt b}){\tt c}\rangle = \langle {\tt b},\diag({\tt a}){\tt c}\rangle = \langle \diag({\tt a}){\tt b},{\tt c}\rangle = \langle {\tt b},\diag({\tt c}) {\tt a} \rangle, \label{eqn:interchanging}
   \end{equation}
   since they all satisfy definition (\ref{eqn:definitionvectorproduct}).
   If relation (\ref{eqn:interchanging}) and the product rule are applied several times, we find
   \begin{align*}
    \frac 12&\naf{}t \langle{\tt v},{\tt rv}\rangle_v
          = \frac 12\naf{}t  \langle{\tt v},  \diag({\sf Interp}_{v\leftarrow c} {\tt rho}) {\tt v} \rangle_v
          = \frac 12\naf{}t   \langle \diag({\tt v}) {\tt v},  {\sf Interp}_{v\leftarrow c} {\tt rho} \rangle_v  \\
         &=            \left\langle \diag({\tt v}) \naf{{\tt v}}t,  {\sf Interp}_{v\leftarrow c} {\tt rho} \right\rangle_v
            + \frac 12 \left\langle \diag({\tt v}) {\tt v}, {\sf Interp}_{v\leftarrow c} \naf {~\tt rho}t \right\rangle_v \\
         &=            \left\langle {\tt v}, \diag\left(\naf{{\tt v}}t\right) {\sf Interp}_{v\leftarrow c} {\tt rho} \right\rangle_v
            \hspace{-0.2cm} +          \left\langle \diag({\tt v}) {\tt v}, {\sf Interp}_{v\leftarrow c} \hspace{-0.1cm} \naf {~\tt rho}t \right\rangle_v
            \hspace{-0.2cm} - \frac 12 \left\langle \diag({\tt v}) {\tt v}, {\sf Interp}_{v\leftarrow c} \hspace{-0.1cm} \naf {~\tt rho}t \right\rangle_v\\
         &=            \left\langle {\tt v}, \naf {~\tt rv}t \right\rangle_v
            - \frac 12 \left\langle {\tt v}, \diag({\tt v}){\sf Interp}_{v\leftarrow c} \naf {~\tt rho}t) \right\rangle_v,
   \end{align*}
   such that
   \[
    \naf{\tt E}t
       =           \left\langle{\tt c1}, \frac{e^\prime_{int}({\tt p})}{R^\prime({\tt p})} \naf{~\tt rho}t \right\rangle_c
         +          \left\langle{\tt v},\naf {~\tt rv}t \right\rangle_v
         - \frac 12 \left\langle{\tt v},\diag({\tt v}) {\sf Interp}_{v\leftarrow c} \naf {~\tt rho}t \right\rangle_v.
   \]
   We can move the factor $e^\prime_{int}({\tt p})/R^\prime({\tt p})$ to the left, and use the continuity and momentum equations (\ref{eqn:Euler}) to replace the time derivatives with spatial derivatives:
   \[ 
    \naf{\tt E}t  =
         -          \left\langle\frac{e^\prime_{int}({\tt p})}{R^\prime({\tt p})}, {\sf DIVr}~{\tt v} \right\rangle_c
         -          \left\langle{\tt v},{\sf GRAD}~{\tt p}\right\rangle_v 
         -          \left\langle{\tt v},{\sf ADVEC}~{\tt v}\right\rangle_v
         + \frac 12 \left\langle{\tt v},\diag({\tt v}) {\sf Interp}_{v\leftarrow c} {\sf DIVr}~{\tt v} \right\rangle_v.
   \]
   The fraction $e^\prime_{int}/R^\prime$ is expanded using the definition of $e_{int}$ for the isentropic Euler equations shown in Table \ref{Tab: mass, momentum and energy}:
   \begin{align*}
       \frac{e^\prime_{int}({\tt p})}{R^\prime({\tt p})}
          &= \frac{1}{R^\prime(p)} \naf{}p \int^p \frac {R(p)-R(q)}{R(q)}~\mbox dq
           = \frac{1}{R^\prime(p)} \naf{}p \left(R(p) \int^p \frac{1}{R(q)}~\mbox dq - \int^p~\mbox dq \right) \\
          &=  \frac{1}{R^\prime(p)} \left(R^\prime(p) \int^p \frac{1}{R(q)}~\mbox dq + R(p) \frac{1}{R(p)} - 1\right) = \int^p \frac{1}{R(q)}~\mbox dq = Q(p).
   \end{align*}
   Furthermore, the chain rule for ${\sf rGRAD}$ is applied (see Table \ref{Tab: properties}), and the two advective terms are factorized:
   \[
    \naf{\tt E}t =
                     - \left\langle Q({\tt p}) , {\sf DIVr}~{\tt v} \right\rangle_c
                     - \left\langle{\tt v},{\sf rGRAD}~Q({\tt p})\right\rangle_v
                     - \left\langle{\tt v},\left( {\sf ADVEC} - \frac 12 \diag\left({\sf Interp}_{v\leftarrow c} {\sf DIVr}~{\tt v}\right)\right) {\tt v} \right\rangle_v.
   \]
   Finally, the symmetry property of ${\sf DIVr}$ and ${\sf rGRAD}$ is applied to the first terms:
   \begin{align*}
     -\left\langle Q({\tt p}) ,{\sf DIVr}~{\tt v} \right\rangle_c - \left\langle{\tt v},{\sf rGRAD}~Q({\tt p})\right\rangle_v
        &= - \left\langle Q({\tt p}) ,{\sf DIVr}~{\tt v} \right\rangle_c
           - \left\langle{\sf rGRAD}^*~{\tt v},Q({\tt p})\right\rangle_v \\
        &= - \left\langle Q({\tt p}) ,{\sf DIVr}~{\tt v} \right\rangle_c
           + \left\langle{\sf DIVr}~{\tt v},Q({\tt p})\right\rangle_v
         = 0,
   \end{align*}
   and also the symmetry property of ${\sf ADVEC}$ (see Table \ref{Tab: properties}):
   \[  \left\langle{\tt v},\left( {\sf ADVEC}
      - \frac 12 \diag\left({\sf Interp}_{v\leftarrow c} {\sf DIVr}~{\tt v}\right)\right) {\tt v} \right\rangle_v
      =
         \left\langle{\tt v},\frac 12 \left( {\sf ADVEC} - {\sf ADVEC}^*\right){\tt v} \right\rangle_v
      = 0, \]
   since it holds that $\langle {\tt x}, ({\sf A} - {\sf A}^*) {\tt x} \rangle = \langle {\tt x}, {\sf A} {\tt x} \rangle - \langle {\tt x}, {\sf A}^* {\tt x} \rangle =  \langle {\tt x}, {\sf A} {\tt x} \rangle - \langle {\sf A}{\tt x}, {\tt x} \rangle = 0$ \correction{for any matrix {\sf A}}{Reviewer 1, Remark 1 and Reviewer 2, Remark 6: Introducing quantities}. This means that \correction{total energy}{Reviewer 2, Remark 7: energy equation}  is indeed conserved: $\naf{\tt E}t  = 0$.

     Now that we have investigated the properties that a discretization should satisfy to be symmetry preserving, we only need to construct such discretizations. This will cover the rest of this paper.

   \begin{table}[ht!]
      \caption{Expressions for the conserved quantities mass, momentum and energy in the three models. Here, $\rho = R(p)$ for the compressible-wave model and isentropic Euler model.}
      \begin{center}
      {\small 
         \begin{tabular}{l|l|l|l|l}
                                     & Linear-wave         & Compressible-wave   & Isentropic Euler  \\
         \hline
                                     &                     &                     &                       \\
         {\bf Mass}                  &                     &                     &                       \\
         ~~~continuous $M$           
                                     & $\displaystyle
                                     \int_V \rho~\mbox dV$ & $\displaystyle
                                                           \int_V \rho~\mbox dV$ & $\displaystyle
                                                                                   \int_V \rho~\mbox dV$ \\[2ex]
                                     &                     &                     &                       \\
         ~~~discrete {\tt M}         
                                     & $\langle{\tt c1},
                                       {\tt rho}\rangle_c$ & $\langle{\tt c1},
                                                             {\tt rho}\rangle_c$ & $\langle{\tt c1},
                                                                                     {\tt rho}\rangle_c$ \\[2ex]
         \hline
                                     &                     &                     &                       \\
         {\bf Momentum}              &                     &                     &                       \\
         ~~~continuous $\vec M$      & $\displaystyle
                                       \rho_0 \int_V \vec v~\mbox dV$ & $\displaystyle
                                                  \rho_0 \int_V \vec v~\mbox dV$ & $\displaystyle
                                                                            \int_V \rho \vec v~\mbox dV$ \\[2ex]
                                     &                     &                     &                       \\
         ~~~discrete $\vec {\tt M}$  & $\rho_0\left(\begin{array}{l}
                                       \langle{\tt c100},{\tt v}\rangle_v \\
                                       \langle{\tt c010},{\tt v}\rangle_v \\
                                       \langle{\tt c001},{\tt v}\rangle_v
                                       \end{array}\right)$ & $\rho_0\left(\begin{array}{l}
                                                             \langle{\tt c100},{\tt v}\rangle_v \\
                                                             \langle{\tt c010},{\tt v}\rangle_v \\
                                                             \langle{\tt c001},{\tt v}\rangle_v
                                                             \end{array}\right)$ & $\left(\begin{array}{l}
                                                                                    \langle{\tt c100},{\tt rv}\rangle_v \\
                                                                                    \langle{\tt c010},{\tt rv}\rangle_v \\
                                                                                    \langle{\tt c001},{\tt rv}\rangle_v
                                                                                    \end{array}\right)$  \\[2ex]
         \hline
                                     &                     &                     &                       \\
         {\bf Energy}                &                     &                     &                       \\
         ~~~continuous $E$           
                                     & $\displaystyle
                                       \int_V  (e_{int} + e_{kin}) ~\mbox dV$ & $\displaystyle
                                            \int_V  (e_{int} + e_{kin}) ~\mbox dV$ &
                                                     $\displaystyle \int_V  (e_{int} + e_{kin}) ~\mbox dV$ \\[3ex]
         ~~~$e_{int}$                & $\displaystyle
                                       \frac {c^2 \rho^2}{2\rho_0}$ & $\displaystyle
                                       \rho_0 \int^p \frac {R(p)-R(q)}{R^2(q)}~\mbox dq $ & $\displaystyle
                                                                \int^p \frac {R(p)-R(q)}{R(q)}~\mbox dq$ \\[2ex]
         ~~~$e_{kin}$                
                                     & $\frac{\rho_0}{2}
                                               |\vec v|^2$ & $\frac{\rho_0}{2}
                                                                     |\vec v|^2$ & $\frac{\rho}{2}
                                                                                             |\vec v|^2$ \\[2ex]
                                     &                     &                     &                       \\
         ~~~discrete ${\tt E}$       
                                     & $\displaystyle \langle
                                       {\tt c1},e_{int}({\tt p})\rangle_c$ & $\displaystyle \langle
                                             {\tt c1},e_{int}({\tt p})\rangle_c$ & $\displaystyle \langle
                                                                     {\tt c1},e_{int}({\tt p})\rangle_c$ \\
                                    &
                                       $\displaystyle  ~~~~~~~~~~               +      \frac{\rho_0}2 \langle{\tt v},{\tt v}\rangle_v$ &
                                       $\displaystyle  ~~~~~~~~~~               +      \frac{\rho_0}2 \langle{\tt v},{\tt v}\rangle_v$ &
                                       $\displaystyle  ~~~~~~~~~~               +       \frac 12 \langle{\tt rv},{\tt v}\rangle_v$
                                       \\[2ex]
                                    &&                      &                     &                       \\
         \end{tabular}
      }\end{center}
   \label{Tab: mass, momentum and energy}
   \end{table}

\section{Numerical approaches}\label{sec:numerical}
   In this section, we describe two different techniques to spatially discretize the models: a Galerkin-type approach (which turns out to be very costly), and a new finite-difference technique. The Galerkin-type approach is only used to introduce some of the concepts that are needed for the finite-difference technique. Once the spatial discretization is found, time integration is needed for a simulation.  We use the ODEPACK solver \verb|lsode| \cite{Hin83}, where it is verified that the time integration is accurate enough such that it has no influence on the results.  Following this approach of highly-accurate time integration allows us to study the spatial discretization separately from the time-integration method, which in actual applications should be symmetry-preserving itself to guarantee the conservation of mass, momentum and/or \correction{total energy}{Reviewer 2, Remark 7: energy equation}.

   \subsection{Galerkin-type approach} \label{sec:Galerkin}
      This section explains an approach that resembles the Galerkin method that is used in finite-element techniques. This technique will be used for the linear-wave model. The method is easy to understand and contains some concepts that help to understand the new finite-difference technique that is explained in the next section.
      
      The traditional Galerkin method samples the solution to obtain a discrete vector of approximations. This discrete approximation is interpolated using a test function, and then the continuous differential operator is applied. The solution is then approximated by the discrete vector for which the inner product of the residual with all test functions equals zero \cite{Vui14R}.

      In order to construct a symmetry-preserving discretization, we follow the same lines. We interpolate the discrete vector of approximations to a continuous function, apply the continuous differential operator, and sample the function back to a discrete vector. If the interpolation and the sampling operator are \emph{mutually adjoint}, and the interpolation of constant functions is exact, then all properties in Table \ref{Tab: properties} (and thereby all conservation laws) are satisfied. The difference between the original Galerkin method and our approach is two-fold: (i) we focus on specific sampling and interpolation operators that are mutually adjoint, and (ii) in the original Galerkin method, the complete residual is sampled (multiplied by a test function), while we keep the original time derivative and only sample the space derivative and right-hand side. This resembles the idea of a lumped mass matrix. Note that the description below is closely related to mimetic discretizations that are defined using \emph{de Rham maps} \cite{Tar99KB}.

      \bigskip
      Let ${\cal J}_c$ be an interpolation operator that maps from the discrete field on the pressure grid to the continuous field, and ${\cal S}_c$ be the sampling operator that produces discrete values from a continuous function \cite{Vet14KG}.  The interpolated fields are written using italic letters and sampled fields with truetype letters, for example, $f := {\cal J}_c~{\tt f}$, and ${\tt g} = {\cal S}_c \ g$. The continuous differential operator ${\cal A}$ is applied to the continuous field obtained from interpolation of the discrete field using ${\cal J}_c$, and the result is mapped back using ${\cal S}_c$, which leads to the discrete operator
      \begin{equation}
        {\sf A} := {\cal S}_c~{\cal A}~{\cal J}_c. \label{eqn:discreteoperator}
      \end{equation}
      The sampling operator and the interpolation operator are called {\em mutually adjoint} if ${\cal S}_c = {\cal J}^*_c$. In the rest of this section, we will see that mutually-adjoint interpolation and sampling, combined with the exact interpolation of the constant functions, is in general enough to obtain the properties of the discrete operators in Table \ref{Tab: properties}. For collocated grids, exact interpolation for constant functions means that ${\cal J}_c ~{\tt c1} = 1$. For staggered grids, we require that the interpolation operator $\vec{\mathcal{J}}_v$ satisfies
      \[
        (\vec{\mathcal{J}}_v {\tt c100})({\vec x}) = (1,0,0), ~~~~
        (\vec{\mathcal{J}}_v {\tt c010})({\vec x}) = (0,1,0), ~~~~
        (\vec{\mathcal{J}}_v {\tt c001})({\vec x}) = (0,0,1).
      \]
      
      \subsubsection{Linear-wave equations: operators {\sf GRAD} and {\sf DIV}}\label{sec:GalerkinScalar}
      In this section, we first define the mutually-adjoint interpolation and sampling operators that are used. Then, we present the symmetry-preserving operators 
      {\sf DIV} and {\sf GRAD} used for the linear-wave equation. Earlier work on the Galerkin-type discretization was presented in \cite{Hof17V}.
      
      For a standard uniform one-dimensional grid with $x_i=i$, standard Lagrange interpolation polynomials, $w$, are used. The corresponding curvilinear interpolation functions are constructed by applying the interpolation function $w$ in each of the three dimensions:
         \begin{equation}
            w_{c,i+m_xj+m_xm_yk}(\vec X(\vec \xi)) =
               w\left( \frac{\xi}{\Delta \xi} - i\right) ~\cdot ~
               w\left( \frac{\eta}{\Delta \eta} - j\right) ~ \cdot ~
               w\left( \frac{\zeta}{\Delta \zeta} - k\right). \label{eqn:itp_function}
         \end{equation}
         The choice of the initial interpolation function $w$ determines the accuracy of the interpolation, as well as the sparseness of the discrete operators. 
         
         To see what the mutual \adjointness\ actually means for the sampling and interpolation operators, they are both written in a more explicit form \cite{Vet14KG}, in terms of the {\em interpolation functions} $w_{c,i}$ ($c$ for cell centers, with point index $i$). As usual, we use
         \begin{subequations}
         \begin{equation}
            ({\cal J}_c{\tt f})(\vec x) = \sum_i {\tt f}_i w_{c,i}(\vec x). \label{eqn:interpolation}
         \end{equation}
         Sampling values will be obtained by calculating the integral of the product of a function and a \emph{sampling function} $s_{c,i}$:
         \begin{equation}
          (\mathcal{S}_c g)_i = \int_V s_{c,i}(\vec x) g(\vec x) \ \mbox dV.
          \label{eqn:sampling}
         \end{equation}
         \end{subequations}
         Mutual adjointness is found using the scalar product on the continuous space (defined by an integral) and the corresponding discrete scalar product on the pressure grid: 
         \begin{align*}
          \langle g, {\cal J}_c {\tt f} \rangle_{\text{continuous}} &= \int_V g(\vec x) (\mathcal{J}_c {\tt f}) (\vec x)\ \mbox dV = \sum_i \int_V g(\vec x) {\tt f}_i w_{c,i}(\vec x)\ \mbox dV = \sum_i {\tt f}_i {\tt dVc}_i \int_V g(\vec x)  \frac{w_{c,i}(\vec x)}{{\tt dVc}_i}\ \mbox dV \\
          &= \sum_i {\tt f}_i {\tt dVc}_i \int_V s_{c,i}(\vec x) g(\vec x)\ \mbox dV = \sum_i {\tt f}_i {\tt dVc}_i (\mathcal{S}_c g)_i = \langle \mathcal{S}_c g, {\tt f} \rangle_c,
         \end{align*}
         such that the mutually-adjoint sampling functions should satisfy
         \begin{equation}
            s_{c,i} = \frac{w_{c,i}(\vec x)}{{\tt dVc}_i} = \frac{w_{c,i}(\vec x)}{\int_V w_{c,i}(\vec x)~\mbox dV}.
         \label{eqn:JS}
         \end{equation}
         Here, we have chosen the integration weights ${\tt dVc}_i$ such that the integral of the sampling function equals 1: in that case, the sampling of a constant field is exact.

         As an example, the standard linear interpolation of a one-dimensional function ${\tt f}$ between ${\tt f}_i$ and ${\tt f}_{i+1}$ is given by
         \[ ({\cal J}_c{\tt f})(x) = {\tt f}_i + \frac{x-x_i}{x_{i+1} - x_i}({\tt f}_{i+1}- {\tt f}_i)
                                   = \frac{x_{i+1} - x}{x_{i+1}-x_i} {\tt f}_i + \frac{x - x_i}{x_{i+1}-x_i} {\tt f}_{i+1}, \quad x_i \leq x \leq x_{i+1}.  \]
         On the whole domain, the interpolated function satisfies equation (\ref{eqn:interpolation}), with 
         \[ w_{c,i}(x) = \left\{  \begin{array}{cl} \frac{x - x_{i-1}}{x_i-x_{i-1}},  &   x_{i-1} \leq x \leq x_i, \\[2ex]
                                                    \frac{x_{i+1}- x}{x_{i+1}-x_i},   &   x_i \leq x \leq x_{i+1}.\end{array}  \right.  \]
         The corresponding sampling operator follows from equation (\ref{eqn:JS}):
         \[ ({\cal S}_c g)_i = \int s_{c,i}(x) g(x) ~\mbox dx =  \int g(x)  \frac{w_{c,i}(x)}{\int w_{c,i}(x)~\mbox dx}  ~\mbox dx
                                                              = \frac{2}{x_{i+1}-x_{i-1}}\int g(x)  w_{c,i}(x)~\mbox dx. \]
         Similarly, higher-order interpolation and sampling operators can be constructed.
         
         Next, we investigate the Galerkin-type approach for vector fields, in order to define the {\sf GRAD} and {\sf DIV} operators for the linear-wave equations. The interpolation of the discrete vector field ${\tt v}$ to the continuous vector field ${\vec v} = \vec{\mathcal{J}}_v {\tt v}$ is very similar to the interpolation for scalar fields (equation (\ref{eqn:interpolation})):
         \begin{equation}
          (\vec{\mathcal{J}}_v {\tt v})({\vec x}) :=
                        {\vec r}_x(\vec x)~({\cal J}_e {\tt vx}) ({\vec x}) +
                        {\vec r}_y(\vec x)~({\cal J}_n {\tt vy}) ({\vec x}) +
                        {\vec r}_z(\vec x)~({\cal J}_t {\tt vz}) ({\vec x}),
         \label{eqn:interpolation_vector_fields}
         \end{equation}
         with
         \begin{equation}
            ({\cal J}_e {\tt vx}) ({\vec x}) := \sum_i {w}_{e,i}({\vec x}) {\tt vx}_i, \quad
            ({\cal J}_n {\tt vy}) ({\vec x}) := \sum_i {w}_{n,i}({\vec x}) {\tt vy}_i, \quad 
            ({\cal J}_t {\tt vz}) ({\vec x}) := \sum_i {w}_{t,i}({\vec x}) {\tt vz}_i, \label{eqn:Jent}
         \end{equation}
         where ${w}_{e,i}$, ${w}_{n,i}$ and ${w}_{t,i}$ are interpolation functions similar to equation (\ref{eqn:itp_function}).

         Sampling of a continuous vector field ${\vec v}$ to its discrete representation ${\tt v} = \mathcal{S}_v {\vec v}$ is similar to the way a continuous scalar field $f$ is sampled to its discrete representation ${\tt f} = {\mathcal S}_c f$ (see equations (\ref{eqn:sampling}) and (\ref{eqn:JS})), and is given by
         \[
          {\tt vx}_i = \int_V \frac{w_{e,i}({\vec x}) {\vec r}_x \cdot {\vec v}_x({\vec x})}{{\tt dVe}_i}\ \mbox dV, \quad {\tt dVe}_i = \int_V w_{e,i}({\vec x})\ \mbox dV,
         \]
         and analogously for ${\tt vy}$ and ${\tt vz}$.

         The Galerkin-type approach requires that exact results are obtained when interpolating the discrete representations {\tt c100}, {\tt c010} and {\tt c001} of the constant vector fields:
         \[
          (\vec{\mathcal{J}}_v {\tt c100})({\vec x}) = (1,0,0), ~~~~
          (\vec{\mathcal{J}}_v {\tt c010})({\vec x}) = (0,1,0), ~~~~
          (\vec{\mathcal{J}}_v {\tt c001})({\vec x}) = (0,0,1).
         \]
         Expanding the interpolation operator in ${\cal J}_v{\tt c100}$ (equation (\ref{eqn:interpolation_vector_fields})) and taking the scalar product with $\vec r_x(\vec x)$, an equation for the interpolation weights $w_{e,i}(\vec x)$ is found. Here, we also use the definition of ${\tt c100}$ (equation (\ref{eqn:c100})), the fact that the grid orientation is orthonormal, and definition (\ref{eqn:Jent}):
         \begin{align}
          {\vec r}_x(\vec x) \cdot (\vec{\mathcal{J}}_v {\tt c100})({\vec x}) &=
                        {\vec r}_x(\vec x) \cdot {\vec r}_x(\vec x)~({\cal J}_e {\tt rxx\_at\_e}) ({\vec x}) +
                        {\vec r}_x(\vec x) \cdot {\vec r}_y(\vec x)~({\cal J}_n {\tt ryx\_at\_n}) ({\vec x})  \nonumber \\  
                      &+
                        {\vec r}_x(\vec x) \cdot {\vec r}_z(\vec x)~({\cal J}_t {\tt rzx\_at\_t}) ({\vec x})
         \nonumber \\ &=
                        ({\cal J}_e {\tt rxx\_at\_e}) ({\vec x})  =  \sum_i w_{e,i}(\vec x)~{\tt rxx\_at\_e}_i.
         \end{align}
         Combining this with the expected result ${\vec r}_x(\vec x) \cdot (\vec{\mathcal{J}}_v {\tt c100})({\vec x})= {\vec r}_x({\vec x}) \cdot (1,0,0)$, and repeating the process for the constant fields {\tt c010} and {\tt c001}, the following three equations are found for the weights $w_{e,i}(\vec x)$:
         \begin{align*}
          \sum_i {w}_{e,i}({\vec x}) {\tt rxx\_at\_e}_i &= {\vec r}_x({\vec x}) \cdot (1,0,0), \\
          \sum_i {w}_{e,i}({\vec x}) {\tt rxy\_at\_e}_i &= {\vec r}_x({\vec x}) \cdot (0,1,0), \\
          \sum_i {w}_{e,i}({\vec x}) {\tt rxz\_at\_e}_i &= {\vec r}_x({\vec x}) \cdot (0,0,1).
         \end{align*}
         This can be written as a matrix-vector equation with three equations and involving one unknown for each point in the $e$-grid:
         \[ \left(
          \begin{array}{c}{\tt rxx\_at\_e}^\top\\
                           {\tt rxy\_at\_e}^\top\\
                           {\tt rxz\_at\_e}^\top \end{array}\right) {w}_{e,:}({\vec x})  = {\vec r}_x({\vec x}).
         \]
         This is an underdetermined system. We solve for small corrections (in the weighted least-squares sense) to the standard Lagrange interpolation polynomials while making sure that the support of the weight function is no larger than that of the Lagrange polynomials. Similar underdetermined systems have to be solved for $w_n$ and $w_t$.

         \bigskip
         The {\sf GRAD} and {\sf DIV} operators in the Galerkin-type approach are given by
         \[
          {\sf GRAD} = \mathcal{S}_v (\nabla \mathcal{J}_c), \quad {\sf DIV} = \mathcal{S}_c (\nabla \cdot \vec{\mathcal{J}}_v ).
         \]
         In practice, the divergence operator in $\vec x_{c,i}$ is computed as follows:
         \[
          {\sf DIV}_i := \sum_{j} ({\sf DIV}_{i,(e,j)} + {\sf DIV}_{i,(n,j)} + {\sf DIV}_{i,(t,j)}),
         \]
         with 
         \[
          {\sf DIV}_{i,(e,j)} := \int_V s_{c,i} (\vec x) \nabla \cdot w_{e,i}(\vec x)\ \mbox{d}V,
         \]
         and similarly for ${\sf DIV}_{i,(n,j)}$ and ${\sf DIV}_{i,(t,j)}$, and where $s_{c,i}$ are the sampling functions of equation (\ref{eqn:JS}).
      
         \subsubsection{Computational costs}\label{sec:computation_Galerkin}
         In the preceding sections, we have presented the Galerkin-type approach to compute discrete operators. In this section, the computational costs are investigated. 
         
         At the start of a simulation, the discrete operators are computed using volume integrals. To guarantee symmetry preservation, these integrals need to be calculated up to machine precision. We achieve high accuracy by using Richardson extrapolation of the composite trapezoidal rule. Although this is a time-consuming process, it is only done once to initialize a simulation. The result is stored in a matrix that corresponds to the discrete operator.
         
         After initialization, the evaluation of the discrete operators is also quite expensive. As an example, we investigate the computational work needed to evaluate the divergence operator, which is related to the number of nonzeros in the {\sf DIV} matrix. The interpolation function $w_{c,\tt ix,iy,iz}$ of order ${\tt N}$ is zero outside the support $\left[ {\tt ix} - {\tt N/2}: {\tt ix}+{\tt N/2}, {\tt iy}-{\tt N/2}:{\tt iy}+{\tt N/2}, {\tt iz}-{\tt N/2}:{\tt iz}+{\tt N/2} \right]$. Similarly, the interpolation function $w_{e,\tt jx+1/2,jy,jz}$ is zero outside the support\\ ${\tt \left[ jx+1/2-N/2: jx+1/2+N/2, jy-N/2:jy+N/2, jz-N/2:jz+N/2 \right]}$.

         Therefore, interpolation functions $w_{c,\tt ix,iy,iz}$ and $w_{e,\tt jx+1/2,jy,jz}$ overlap for all indices {\tt jx,jy,jz} for which
               \begin{align*}
                   \tt jx &\in \tt ix-N  :ix+N-1,  \\
                   \tt jy &\in \tt iy-N+1:iy+N-1,  \\
                   \tt jz &\in \tt iz-N+1:iz+N-1, 
               \end{align*}
         which is ${\tt 2*N*(2*N-1)*(2*N-1)}$ functions. The divergence matrix has three times as many nonzeros (for $v_x$, $v_y$ and $v_z$), which is
         ${\tt 6*N*(2*N-1)*(2*N-1)}$ ${\tt \approx 24*N^3}$ nonzeros per matrix row. In 2D, we find ${\tt 4*N*(2*N-1) \approx 8*N^2}$ nonzeros per matrix row, and in 1D ${\tt 2*N}$. 

         A cheaper approach to arrive at symmetry-preserving discretizations is presented in the next section.

   \subsection{Finite-difference approach}\label{sec:FD}
      The Galerkin-type approach described in the previous section has certain disadvantages. For example, the calculation of the discrete operators requires highly-accurate integration, and the resulting matrices have many nonzeros.  Therefore, this method leads to relatively long calculation times. Furthermore, the application of the Galerkin-type approach to difficult operators such as the advection operator is not straightforward.
      
      A more efficient and flexible approach is found using a finite-difference strategy. In this section, we explain the finite-difference technique for the different operators needed in the linear-wave, compressible-wave and isentropic Euler models.

      \subsubsection{Linear-wave equations: operators {\sf DIV} and {\sf GRAD}}\label{Sec: div at point i}

         \correction{A higher-order finite-difference approximation {\sf DIV} for the divergence is constructed from an exact formulation (\ref{eqn:exactdiv}) of the divergence in a cell center.  Subsequently, standard finite-difference stencils are applied to this formulation, resulting in a higher-order, conservative approximation of the divergence. The resulting discrete operator {\sf DIV} combines a series of steps such as interpolation and differentiation with respect to $\vec \xi$. All of these steps have continuous equivalents, and therefore its negative adjoint ${\sf GRAD}:=-{\sf DIV}^*$ is an accurate approximation of the gradient. The approach is very similar to the Richardson extrapolation scheme used in \cite{Kok06} and \cite{Kok09}, but the current approach leads to smaller stencils, especially for high orders of accuracy.}{Reviewer 1, Remark 5: Literature} 

         \bigskip
   \correction{
       The resulting finite-difference discretizations are locally conservative for the mass balance, because the change in a control volume can be written as the net effect of fluxes, each of which describes the transport between two 'neighbor' grid points.
       In rectilinear, orthogonal, possibly stretched, grids the momentum equation is locally conservative in the same sense, with momentum control volumes drawn around the cell faces of the mass control volumes.
       In curvilinear or nonorthogonal grids, the momentum equation cannot be seen as locally conservative, because every sampled velocity contains momentum in a different direction, making it impossible to transfer momentum from one grid point/control volume to another.  
       For similar reasons, the energy balance cannot be seen as locally conservative, because internal energy is sampled in the mass control volumes, and kinetic energy in the cell faces. 
       }{ Reviewer 1, Remark 6: locally conservative} 

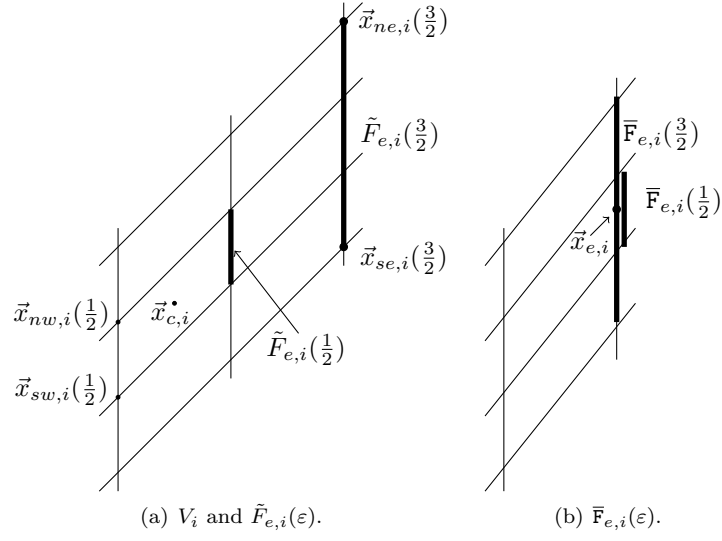
\begin{figure}[ht!]
        \centering
        \subfigure[$V_i$ and $\tilde F_{e,i}(\eps)$.]{
        \begin{tikzpicture}[scale = 0.5]
       \draw (0,-2.5)--(0,4.5);
       \draw (3,0.5)--(3,7.5);
       \draw[line width=2pt] (3,3)--(3,5);
       \draw (-0.5,-2.5)--(6.5,4.5);
       \draw (-0.5,-0.5)--(6.5,6.5);
       \draw (-0.5,3.5)--(6.5,10.5);
       \draw (6,3.5)--(6,10.5);
       \draw[line width=2pt] (6,4)--(6,10);
       \draw (-0.5,1.5) -- (6.5,8.5);
       \draw[fill=black] (1.5,2.5) circle (0.5mm);
       \node at (1.4,2.2) {$\vec x_{c,i}$};
       \draw[<-] (3.1,3.9) -- (4.8,1.7);
       \node at (5,1.3) {$\tilde F_{e,i}(\frac12)$};
       \node at (7.5,7) {$\tilde F_{e,i}(\frac32)$};
       \draw[fill=black] (0,2) circle (0.5mm);
       \node at (-1.5,2.2) {$\vec x_{nw,i}(\frac12)$};
       \draw[fill=black] (0,0) circle (0.5mm);
       \node at (-1.5,0.2) {$\vec x_{sw,i}(\frac12)$};
       \draw[fill=black] (6,10) circle (1mm);
       \node at (7.6,10) {$\vec x_{ne,i}(\frac32)$};
       \draw[fill=black] (6,4) circle (1mm);
       \node at (7.6,3.7) {$\vec x_{se,i}(\frac32)$};
       \end{tikzpicture}\label{fig:tildeF}
       }
        \subfigure[$\overline {\tt F}_{e,i}(\eps)$.]{
        \begin{tikzpicture}[scale = 0.5]
       \draw (0,-2.5)--(0,4.5);
       \draw (3,1)--(3,8.5);
       \draw (-0.5,-2.5)--(3.5,2.5);
       \draw (-0.5,-0.5)--(3.5,4.5);
       \draw (-0.5,3.5) --(3.5,8.5);
       \draw (-0.5,1.5) --(3.5,6.5);
       \draw[line width=2pt] (3.2,4)--(3.2,6);
       \draw[line width=2pt] (3,2)--(3,8);
       \draw[fill=black] (3,5) circle (1mm);
       \draw[<-] (2.8,4.8)--(2.3,4.3);
       \node at (2.3,4.1) {$\vec x_{e,i}$};
       \node at (4.8,5.2) {$\overline {\tt F}_{e,i}(\frac12)$};
       \node at (4.2,7) {$\overline {\tt F}_{e,i}(\frac32)$};
       \end{tikzpicture}\label{fig:overlineF}
       }
       \caption{Controle volumes and approximate fluxes (over bold faces). For $\eps = 1/2$ (in Figure \ref{fig:tildeF}), the control volume equals the original grid cell, and $\tilde F_{e,i}(\eps) = \overline {\tt F}_{e,i}(\eps)$.}\label{fig:control_volume}
      \end{figure}

         \bigskip
         For the derivation of the finite-difference technique, a control volume $V_i$ is used, centered around $\vec x_{c,i}$.
         The control volume is an interval (1D), quadrilateral (2D) or hexahedron (3D).
         The method is described in 2D for simplicity, but it works the same way in 1D and 3D.
         In 2D, the control volume $V_i$ is given by
         \[  V_i(\eps) := \vec x_{ne,i}(\eps)\vec x_{nw,i}(\eps)\vec x_{sw,i}(\eps)\vec x_{se,i}(\eps), \]
         with vertices given by
         \begin{align*}
            \vec x_{ne,i}(\eps) := \vec X\left(\vec \xi_{c,i} + \eps \mbox{\scriptsize{$\left(\!\!\begin{array}{c}\Delta \xi\\\Delta \eta \end{array}\!\!\right)$}}\right),& \quad
            \vec x_{nw,i}(\eps) := \vec X\left(\vec \xi_{c,i} + \eps \mbox{\scriptsize{$\left(\!\!\begin{array}{r}-\Delta \xi\\\Delta \eta \end{array}\!\!\right)$}}\right),
          \\
            \vec x_{sw,i}(\eps) := \vec X\left(\vec \xi_{c,i} + \eps \mbox{\scriptsize{$\left(\!\!\begin{array}{r}-\Delta \xi\\-\Delta \eta \end{array}\!\!\right)$}}\right),& \quad
            \vec x_{se,i}(\eps) := \vec X\left(\vec \xi_{c,i} + \eps \mbox{\scriptsize{$\left(\!\!\begin{array}{r}\Delta \xi\\-\Delta \eta \end{array}\!\!\right)$}}\right),
         \end{align*}
         see Figure \ref{fig:tildeF}.
         
         The scaled flux of the vector field $\vec v$ out of the east cell face is denoted $F_{e,i}$ and given by
         \[ F_{e,i}(\eps) := \frac 1\eps\int_{\vec x_{se,i}(\eps)}^{\vec x_{ne,i}(\eps)} \vec v(\vec x) \cdot \vec n ~\mbox dS, \]
         and the scaled fluxes $F_{w,i}=-F_{e,i}(-\eps)$, $F_{n,i}$ and $F_{s,i}=-F_{n,i}(-\eps)$ out of the west, north and south cell faces analogously.

         Gauss's theorem equates the integrated divergence to the net outflux:
         \[
            \int_{V_i(\eps)} \hspace{-0.6cm}\nabla \cdot \vec v(\vec x) ~\mbox dV \hspace{-0.1cm}= \eps(F_{e,i}(\eps)+F_{w,i}(\eps)+F_{n,i}(\eps)+F_{s,i}(\eps)) 
                                                                   = \eps(F_{e,i}(\eps)-F_{e,i}(-\eps)+F_{n,i}(\eps)-F_{n,i}(-\eps)).
         \]
         For very small $\eps$, the volume integral may be approximated by the midpoint rule (note the occurrence of the Jacobian determinant of the transformation $\vec{X}$):
         \begin{align*}
            4 \eps^2 \Delta \xi \Delta \eta \left| \paf{\vec X}{\vec \xi}\right| \nabla \cdot \vec v(\vec x_{c,i}) + {\cal O}(\eps^4) &=
         \int_{V_i(\eps)} \nabla \cdot \vec v(\vec x) ~\mbox dV \\
         &=\eps(F_{e,i}(\eps)-F_{e,i}(-\eps)+F_{n,i}(\eps)-F_{n,i}(-\eps)),
         \end{align*}
         and this leads to the following {\em exact} representation of the divergence $\nabla \cdot \vec v$ (when taking $\eps \rightarrow 0$):
         \[
         \nabla \cdot \vec v(\vec x_{c,i})  =
            \frac {F^\prime_{e,i}(0) + F^\prime_{n,i}(0)}
         {2\Delta \xi\Delta \eta \left| \paf{\vec X}{\vec \xi}\right|} =
            \frac {F^\prime_{e,i}(0) + F^\prime_{n,i}(0)} {2{\tt dVc}_i},
         \]
         because the integration weights will be chosen equal to
         \[
            {\tt dVc}_i := \Delta \xi\Delta \eta\left|\paf{\vec X}{\vec \xi}(\vec x_{c,i})\right|.
         \]

         Note that the scaled fluxes $F_{e,i}$ and $F_{n,i}$ are given in terms of boundary integrals. They can be replaced by approximations
         in which only velocity values $\vec v$ occur, and no integration is needed.
         To do this, the cell-face centers $\vec x_{e,i}$ and $\vec x_{n,i}$ are defined by
         \[
         \vec x_{e,i}(\eps) := \vec X(\vec \xi_{e,i}(\eps)) := \vec X\left(\vec \xi_{c,i}+\eps \mbox{\scriptsize{$\left(\!\!\begin{array}{c} \Delta \xi\\0\end{array}\!\!\right)$}}\right), \quad
         \vec x_{n,i}(\eps) := \vec X(\vec \xi_{n,i}(\eps)) := \vec X\left(\vec \xi_{c,i}+\eps \mbox{\scriptsize{$\left(\!\!\begin{array}{c} 0\\\Delta \eta\end{array}\!\!\right)$}}\right).
         \]
         The flux $F_{e,i}$ is approximated by the approximate flux $\tilde F_{e,i}$, obtained using the midpoint rule:
         \begin{align}
            \tilde F_{e,i}(\eps) :=& \frac 1\eps\int_{\vec x_{se,i}(\eps)}^{\vec x_{ne,i}(\eps)} \vec n  ~\mbox dS ~\cdot~ \vec v(\vec x_{e,i}(\eps)) \nonumber \\
                                 =& \frac 1\eps |\vec x_{ne,i} - \vec x_{se,i}| \frac{1}{|\vec x_{ne,i} - \vec x_{se,i}|}
                                    \begin{pmatrix}
                                       y_{ne,i} - y_{se,i} \\ -(x_{ne,i} - x_{se,i}) \\
                                    \end{pmatrix}
                                    ~\cdot~ \vec v(\vec x_{e,i}(\eps)) \nonumber \\
                                 =&  -\frac {(\vec x_{ne,i}(\eps)-\vec x_{se,i}(\eps))^\perp}\eps \cdot  \vec v(\vec x_{e,i}(\eps)),\label{eqn:Fe}
         \end{align}
         where the definition of the outward-directed normal on the east cell face is used, and the {\em perp operator} $\perp$ is defined as $(x,y)^\perp := (-y,x)$. The location of $\tilde F_{e,i}(\eps)$ for two different values of $\eps$ is visualized in Figure \ref{fig:tildeF}.
         
         Similarly, we define
         \begin{equation}
            \tilde F_{n,i}(\eps) := \frac {(\vec x_{ne,i}(\eps)-\vec x_{nw,i}(\eps))^\perp}\eps \cdot v(\vec x_{n,i}(\eps)). \label{eqn:Fn}
         \end{equation}
         Since the midpoint rule is at least second-order accurate, the approximate fluxes and their derivatives are exact at $\eps = 0$.
         Hence, an {\em exact} formulation for the divergence $\nabla \cdot \vec v$ is given by
         \begin{equation}
         \nabla \cdot \vec v(\vec x_{c,i})  =
         \frac{\tilde F^\prime_{e,i}(0) + \tilde F^\prime_{n,i}(0)}{2{\tt dVc}_i}, \label{eqn:exactdiv}
         \end{equation}
         which indeed does not include any integrals.
         
         \bigskip
         The discrete divergence ${\sf DIV}$ is now found by executing the following steps: 
         \begin{enumerate}
          \item Interpolation of the velocities.
          
           To compute $\tilde{F}_{e,i}(\eps)$ (equation (\ref{eqn:Fe})), we need to evaluate the discrete value of 
           \[
            v(\vec x_{e,i}(\eps)) = \vec r_x(\vec x_{e,i}(\eps)) v_x(\vec x_{e,i}(\eps)) + \vec r_y(\vec x_{e,i}(\eps)) v_y(\vec x_{e,i}(\eps)),
           \]
           see equation (\ref{eq: v from components}). Since $\vec x_{e,i}(\eps)$ is a point in the $e$-grid, the discrete $y$-component of $\vec v$ (the sampled velocity components {\tt vy}) is not available yet. The component ${\tt vy}$ is first interpolated to the pressure grid (using the {\em destaggering matrix} {\sf N2C}) and then to the $e$-grid (with the {\em staggering matrix} {\sf C2E}) using an interpolation procedure similar to the one in Section \ref{sec:Galerkin}.  The result will be stored in the vector ${\tt vy\_at\_e} = {\sf C2E}~{\sf N2C}~{\tt vy} =: {\sf N2E}~{\tt vy}$.

           Similarly, the components {\tt vx} are interpolated to ${\tt vx\_at\_n}={\sf C2N}~{\sf E2C}~{\tt vx} =: {\sf E2N}~{\tt vx}$, so that the complete velocity vector is available at $e$- and $n$-points.  
           
           In matrix-vector notation, this interpolation step can be written as
           \begin{equation}
            \begin{pmatrix}
             {\tt vx} \\
             {\tt vy\_at\_e} \\
             {\tt vx\_at\_n} \\
             {\tt vy}
            \end{pmatrix}
             =
            \begin{pmatrix}
             {\sf I} & {\sf 0} \\
             {\sf 0} & {\sf N2E} \\
             {\sf E2N} & {\sf 0} \\
             {\sf 0} & {\sf I}
            \end{pmatrix}
            \begin{pmatrix}
             {\tt vx} \\
             {\tt vy}
            \end{pmatrix}. \label{eqn:interpolation_v}
           \end{equation}
           \item Computation of the fluxes.
           
             All the approximate fluxes use the velocity vector in one point of the staggered grid: every $\tilde F_{e,i}$ uses the velocity vector at $e$-points, and $\tilde F_{n,i}$ at $n$-points. Therefore, the flux values are located at these grid points, for instance, the flux at the $e$-point $\vec x_{e,i}$ is
             \begin{align*}
             \overline {\tt F}_{e,i}(\eps) &:= 
             -\frac {
                  \left(\vec X\left(\vec \xi_{e,i} + \eps \mbox{\scriptsize{$\left(\!\! \begin{array}{c} 0 \\ \Delta \eta \end{array}\!\!\right)$}}\right) -
                        \vec X\left(\vec \xi_{e,i} - \eps \mbox{\scriptsize{$\left(\!\! \begin{array}{c} 0 \\ \Delta \eta \end{array}\!\!\right)$}}\right)\right)^\perp
             }{\eps} \cdot
                         \left(
                                 \vec r_x(\vec x_{e,i}) {\tt vx}_i        +
                                 \vec r_y(\vec x_{e,i}) {\tt vy\_at\_e}_i
                         \right)
              \\ &=:
                 {\tt Nx}_{e,i}(\eps)~{\tt vx}_i  +  {\tt Ny}_{e,i}(\eps) ~ {\tt vy\_at\_e}_i,
             \end{align*}
             and similar for $\overline {\tt F}_{n,k}$ for the $n$-point $\vec x_{n,k}$. In matrix-vector form this equals
             \begin{equation}
              \begin{pmatrix}
               \overline {\tt F}_e(\eps) \\
               \overline {\tt F}_n(\eps)
              \end{pmatrix}
              =
              \begin{pmatrix}
               \diag({\tt Nx}_e(\eps)) & \diag({\tt Ny}_e(\eps)) & {\sf 0} & {\sf 0} \\
               {\sf 0} & {\sf 0} & \diag({\tt Nx}_n(\eps)) & \diag({\tt Ny}_n(\eps))
              \end{pmatrix}
              \begin{pmatrix}
                {\tt vx} \\
                {\tt vy\_at\_e} \\
                {\tt vx\_at\_n} \\
                {\tt vy}
              \end{pmatrix}. \label{eqn:matrixFen}
             \end{equation}
             Note that the use of the interpolated velocity vectors is not the only difference between $\tilde{F}_{e,i}$ and $\overline {\tt F}_{e,i}$. Also the $x$-coordinate of 
             \[
              \vec x_{ne,i}(\eps)-\vec x_{se,i}(\eps) =  \vec X\left(\vec \xi_{c,i} + \eps \mbox{\scriptsize{$\left(\!\!\begin{array}{c}\Delta \xi\\\Delta \eta \end{array}\!\!\right)$}}\right) - \vec X\left(\vec \xi_{c,i} + \eps \mbox{\scriptsize{$\left(\!\!\begin{array}{r}\Delta \xi\\-\Delta \eta \end{array}\!\!\right)$}}\right),
             \]
             used in $\tilde{F}_{e,i}$ (at distance $\eps \Delta \xi$, see Figure \ref{fig:tildeF}), differs from the $x$-coordinate of
             \[
              \vec X\left(\vec \xi_{e,i} + \eps \mbox{\scriptsize{$\left(\!\! \begin{array}{c} 0 \\ \Delta \eta \end{array}\!\!\right)$}}\right) -
                        \vec X\left(\vec \xi_{e,i} - \eps \mbox{\scriptsize{$\left(\!\! \begin{array}{c} 0 \\ \Delta \eta \end{array}\!\!\right)$}}\right),
             \]
             used in $\overline {\tt F}_{e,i}$ (fixed at cell faces, see Figure \ref{fig:overlineF}). For constant vector fields, the interpolation is chosen to be exact, and we find the following relation:
             \[
              \tilde{F}_{e,i}(\eps) = \overline{\tt F}_{e,i+\eps-1/2}(\eps), \quad \overline{\tt F}_{e,i}(\eps) = \tilde{F}_{e,i-\eps+1/2}(\eps), \quad \eps\in \left\{ \cdots, -3/2, -1/2, 1/2, 3/2, \cdots \right\}
             \]
             For non-constant vectors, this only approximately holds. Similarly, the following relation for constant vector fields holds: 
             \[
               \tilde{F}_{n,i}(\eps) = \overline{\tt F}_{n,i+(\eps-1/2)m_x}(\eps), \quad \overline{\tt F}_{n,i}(\eps) = \tilde{F}_{e,i-(\eps+1/2)m_x}(\eps).
             \]
             
           \item Finite-difference step.
           
              The flux vectors $\overline{\tt F}_e = (\overline{\tt F}_{e,j})_j$ and $\overline{\tt F}_n = (\overline{\tt F}_{n,k})_k$ are used to calculate the discrete operator {\sf DIV} corresponding to equation (\ref{eqn:exactdiv}). Standard differentiation stencils are applied to calculate {\sf DIV} (which lives on the pressure points). Therefore, we use the matrices ${\sf DIFFX}$ and ${\sf DIFFY}$, that make sure that central differences are used. In internal pressure points $\vec x_{c,i}$, row $i$ of {\sf DIFFX} contains a $-1$ at position $(i,i-(\eps+1/2))$ and a $1$ at $(i,i+(\eps-1/2))$, and {\sf DIFFY} is defined in a similar way, so that the $i$-th entry is given by
             \begin{align*}
             ({\sf DIFFX(\eps)}~\overline {\tt F}_e(\eps))_i &= \overline {\tt F}_{e,i+\eps-1/2} -  \overline {\tt F}_{e,i-\eps-1/2}, \\
             ({\sf DIFFY(\eps)}~\overline {\tt F}_n(\eps))_i &= \overline {\tt F}_{n,i+(\eps-1/2)m_x} -  \overline {\tt F}_{n,i-(\eps+1/2)m_x}.
             \end{align*}
             The divergence (on the whole domain) can be computed as follows
             \begin{equation}
             {\sf DIV}~{\tt v} =
             \frac 12 \diag({\tt dVc})^{-1}
             \sum_{\eps=1/2,~3/2,~\ldots}  \left(
                                                      {\sf DIFFX}(\eps) \overline {\tt F}_e(\eps) +
                                                      {\sf DIFFY}(\eps) \overline {\tt F}_n(\eps) \right) \alpha(\eps),
             \label{eq: discrete DIV}
             \end{equation}
             where $\alpha$ contains the coefficients of the differentiation stencil. For example, for second-order differentiation, $\alpha(1/2) = 1$, and for fourth order, we use \cite{Kin90KS, Liu09S}
             \begin{eqnarray}
                \alpha\left(1/2\right) = \frac 98 &,&
                \alpha\left(3/2\right) = -\frac 1{24}.
             \end{eqnarray}
         \end{enumerate}   
        
         We need to show that this choice for the divergence operator indeed satisfies the null space properties (Table \ref{Tab: properties}). The left null space property holds if ${\sf DIV}^* {\tt c1} = 0$. Indeed, the integral of a divergence is zero:
               \[
                   \left\langle {\tt c1}, {\sf DIV}~{\tt v}\right\rangle_c =
                    {\tt c1}^\top \diag({\tt dVc}) {\sf DIV}~{\tt v} =
               \frac 12
               {\tt c1}^\top \hspace{-0.6cm}
               \sum_{\eps=1/2,~3/2,~\ldots} \hspace{-0.5cm} \left( 
                                                            {\sf DIFFX}(\eps) \overline {\tt F}_e(\eps) +
                                                            {\sf DIFFY}(\eps) \overline {\tt F}_n(\eps) \right) \alpha(\eps)=0,
               \]
               because ${\tt c1}^\top {\sf DIFFX}(\eps) = {\tt c1}^\top {\sf DIFFY}(\eps) = 0$. Next, we prove that the divergence of a constant vector field is zero. The interpolations used to calculate {\tt vx\_at\_n} and {\tt vy\_at\_e} are constructed such that they are exact for constant vector fields. Therefore, for a
               discrete vector field {\tt v} representing a constant vector field ($\vec v(\vec x) = \vec c$), a zero divergence is found, because at grid point $\vec x_{c,i}$, we have (note the exact relation between $\overline{\tt F}$ and $\tilde{F}$)
               \begin{align*}
               \left(
                     {\sf DIFFX}(\eps) \overline {\tt F}_e(\eps) +
                     {\sf DIFFY}(\eps) \overline {\tt F}_n(\eps) \right)_i =& \
               \tilde F_{e,i}(\eps) - \tilde F_{e,i}(-\eps) +
               \tilde F_{n,i}(\eps) - \tilde F_{n,i}(-\eps)
               \nonumber \\  =&
                  -\frac{(\vec x_{ne,i}(\eps) - \vec x_{se,i}(\eps))^\perp} \eps \cdot \vec c
                  +\frac{(\vec x_{nw,i}(\eps) - \vec x_{sw,i}(\eps))^\perp} \eps \cdot \vec c
               \nonumber \\ &+
                   \frac{(\vec x_{ne,i}(\eps) - \vec x_{nw,i}(\eps))^\perp} \eps \cdot \vec c
                  -\frac{(\vec x_{se,i}(\eps) - \vec x_{sw,i}(\eps))^\perp} \eps \cdot \vec c
               = 0.
               \end{align*}

         \bigskip
         The gradient operator ${\sf GRAD}$ is obtained from the divergence's adjoint and is given by
         \begin{equation}
             {\sf GRAD} := -{\sf DIV}^*.
         \label{eq: GRAD}
         \end{equation}
         Using the notations ${\tt Nx}_e$, ${\tt Nx}_n$, ${\tt Ny}_e$, ${\tt Ny}_n$, {\sf E2N} and {\sf N2E}, introduced in this section, and the definition of the adjoint,
         the gradient operator {\sf GRAD} is given by
         \[
             {\sf GRAD} =
                   \sum_{\eps=1/2,~3/2,~\ldots} \frac{\alpha(\eps)}2{\sf GRAD1}(\eps)
                  \left(  \begin{array}{c}
                       -{\sf DIFFX}^\top(\eps) \\
                       -{\sf DIFFY}^\top(\eps)
                  \end{array} \right),
         \]
         where the operator {\sf GRAD1} is given by
         \[
             {\sf GRAD1}(\eps) =
                   \diag
                   \left( \begin{array}{cccc}
                       {\tt dVe}\\
                       {\tt dVn}
                   \end{array}\right)^{-1}
                   \left( \begin{array}{cccc}
                       {\sf I} & {\sf 0} &       {\sf E2N}^\top & {\sf 0} \\ \\
                       {\sf 0} & {\sf N2E}^\top & {\sf 0} & {\sf I}
                   \end{array}\right)
                                        \left( \begin{array}{cccc}
                                               \diag({\tt Nx}_{e}(\eps)) & {\sf 0} \\
                                               \diag({\tt Ny}_e(\eps)) & {\sf 0} \\
                                               {\sf 0} & \diag({\tt Nx}_{n}(\eps)) \\
                                               {\sf 0} & \diag({\tt Ny}_n(\eps))
                                        \end{array}\right).
         \]
         In general, one cannot expect the (discrete) adjoint ${\sf GRAD} = -{\sf DIV}^*$ of an accurate discretization of the divergence to be an accurate discretization of the gradient.  However, the divergence {\sf DIV} can be seen as an operator in computational space, discretized on a uniform $(\xi,\eta)$-grid.  On uniform grids, the adjoint of an interpolation operator is itself an accurate interpolation (in the opposite direction).  Similarly, the adjoint of a discrete derivative operator is itself an accurate discrete derivative operator on uniform grids. Therefore, an accurate discretization of the continuous gradient operator is found if adjoint interpolations and adjoint derivative operators are combined to a (discrete) {\sf GRAD} operator.

         \bigskip
         The finite-difference approach is computationally more efficient than the Galerkin-type approach. Similar to Section \ref{sec:computation_Galerkin}, we investigate the computational work needed to evaluate the divergence operator. Here, we assume interpolation and differentiation of order {\tt N}. We start with the interpolation of $v_x$, $v_y$ and $v_z$ to the pressure grid points. Each interpolation involves {\tt N} grid points, so this step is a matrix-vector multiplication with {\tt 3 * N} nonzeros per grid point.
         Next, $v_y$ and $v_z$ are interpolated to the points of the $e$-grid,
         $v_x$ and $v_z$ to the points of the $n$-grid,  and
         $v_x$ and $v_y$ to the points of the $t$-grid. This is three matrix-vector multiplications, each involving ${\tt 2 * N}$ nonzeros. Finally, the evaluation of the divergence involves an $x$-, $y$- and $z$-derivative, each of which involves ${\tt N}$ nonzeros. Therefore, the evaluation of the divergence requires ${\tt 3 * N + 3 * 2 * N + 3 * N }= {\tt 12 * N}$ nonzeros per grid point. In 2D, the calculation involves ${\tt 2 * N} + {\tt 2 * N + 2 * N} = {\tt 6 * N}$ nonzeros per grid point, and in 1D {\tt N} nonzeros. Indeed, the calculation of the divergence operator is much cheaper than when using a Galerkin-type approach. This observation generally holds for all the operators that are designed in this paper.

      \subsubsection{Compressible-wave equations: operators $\tilde{\sf r}{\sf GRAD}$ and ${\sf DIV}\tilde{\sf r}$}\label{sec:rGRAD}
         The discretization of the compressible-wave equations requires the operator $\tilde{\sf r}{\sf GRAD}$, which is an approximation of the operator $(\tilde{\sf r}{\sf GRAD}~{\tt f})_i\approx (\rho \nabla f)(\vec x_{c,i})$. The operator $\tilde{\sf r}{\sf GRAD}$ is obtained by applying the same interpolations used in the {\sf GRAD}-operator,
         \begin{align*}
            \tilde{\sf r}{\sf GRAD} &:=
                   \sum_{\eps=1/2,~3/2,~\ldots} \frac{\alpha(\eps)}2{\sf GRAD1}(\eps) ~
                   \diag
                   \left( \begin{array}{cccc}
                       \widetilde{\tt rho}_e(\eps)\\
                       \widetilde{\tt rho}_n(\eps)
                   \end{array}\right)
                  \left(  \begin{array}{c}
                       -{\sf DIFFX}^\top(\eps) \\
                       -{\sf DIFFY}^\top(\eps)
                  \end{array} \right), \\
                  {\sf DIV}\tilde{\sf r} &:= - \tilde{\sf r}{\sf GRAD}^*,
         \end{align*}
         where the intermediate densities $\widetilde{\tt rho}_e$ and $\widetilde{\tt rho}_n$  are given by
         \[
           \diag \left( \begin{array}{cccc}
                \widetilde{\tt rho}_e(\eps)\\
                \widetilde{\tt rho}_n(\eps)
            \end{array}\right) :=
            \diag    \left(  \begin{array}{c}
                       {\sf DIFFX}^\top(\eps)~Q({\tt p})\\
                       {\sf DIFFY}^\top(\eps)~Q({\tt p})
                  \end{array} \right)
            \diag\left(
                  \begin{array}{c}
                       {\sf DIFFX}^\top(\eps)~S({\tt p}) \\
                       {\sf DIFFY}^\top(\eps)~S({\tt p})
                  \end{array} \right)^{-1}
              .
         \]
         Here, the functions $Q$ and $S$ that were introduced in Table \ref{Tab: properties} are applied. Obviously, the chain rule for $\tilde{\sf r}{\sf GRAD}$ holds: $\tilde{\sf r}{\sf GRAD} \ S({\tt p}) = {\sf GRAD}\ Q({\tt p})$.
         Of course, special care is necessary for the numerically-stable calculation of the intermediate densities, especially when elements of the diagonal matrix $\diag({\sf DIFF}~S({\tt p}))$ are very small.
         
   \subsection{Isentropic Euler equations: operators {\sf rGRAD} and {\sf DIVr}}\label{sec:rGRAD and DIVr}
      The techniques presented in the preceding sections are almost sufficient for a symmetry-preserving discretization of the isentropic Euler equations.  Only the advection operator is not yet available. The advection operator is very similar to the divergence operator, on which it will be based.  Here, we do not use the operators ${\sf DIV}\tilde{\sf r}$ and $\tilde{\sf r}{\sf GRAD}$, but the very similar operators ${\sf DIVr}$ and ${\sf rGRAD}$, given by
      \begin{align*}
         {\sf rGRAD} &:=
                \sum_{\eps=1/2,~3/2,~\ldots} \frac{\alpha(\eps)}2{\sf GRAD1}(\eps) ~
                   \diag
                   \left( \begin{array}{cccc}
                       {\tt rho}_e(\eps)\\
                       {\tt rho}_n(\eps)
                   \end{array}\right)
               \left(  \begin{array}{c}
                    -{\sf DIFFX}^\top(\eps) \\
                    -{\sf DIFFY}^\top(\eps)
               \end{array} \right), \nonumber \\
               {\sf DIVr} &:= -{\sf rGRAD}^*,
      \end{align*}
         where the intermediate densities ${\tt rho}_e$ and ${\tt rho}_n$  are given by
         \[
           \diag \left( \begin{array}{cccc}
                {\tt rho}_e(\eps)\\
                {\tt rho}_n(\eps)
            \end{array}\right) :=
            \diag     \left(  \begin{array}{c}
                       {\sf DIFFX}^\top(\eps)~{\tt p}\\
                       {\sf DIFFY}^\top(\eps)~{\tt p}
                  \end{array} \right)
            \diag\left(
                  \begin{array}{c}
                       {\sf DIFFX}^\top(\eps)Q({\tt p}) \\
                       {\sf DIFFY}^\top(\eps)Q({\tt p})
                  \end{array} \right)^{-1}
             .
         \]

   \subsection{Isentropic Euler equations: operator {\sf ADVEC}}\label{sec:Advection}

         \correction{A symmetry-preserving approximation for the advection term, which also preserves momentum, is constructed using a 'splitting' strategy very similar to the one presented in \cite{Win18WGW}, \cite{Win17WGK} and \cite{Hof12V}, and classical results of Tadmor \cite{Tad03}. For the construction of the advection operator {\sf ADVEC}, we start with a discrete advection operator ${\sf ADVEC}_s$ that works on a scalar field ${\tt f}_v$. Unlike a 'normal' scalar field {\tt f}, which is sampled at the cell centers, the scalar field  ${\tt f}_v = ({\tt f}_e^\top, {\tt f}_n^\top)^\top$ is sampled at the cell faces. The scalar advection operator ${\sf ADVEC}_s$ approximates the continuous advection operator: $({\sf ADVEC}_s {\tt f}_v)_{c,i} \approx (\nabla \cdot \rho \vec v f)(\vec x_{c,i})$. This operator ${\sf ADVEC}_s$ has the following properties compared to {\sf ADVEC}:
         \begin{itemize}
            \item {\bf Input field of the operator}:

                 ${\sf ADVEC}_s$ works on scalar fields sampled at the cell faces, and not, like {\sf ADVEC}, on a staggered representation of a vector field;

            \item {\bf Output field of the operator}:
            
                 The result ${\sf ADVEC}_s~{\tt f}_v$ lives on cell centers, and is not, like the result ${\sf ADVEC}~{\tt v}$, a staggered representation of a vector field;

            \item {\bf Applying the operator to a constant field}:
            
                When applied to a constant field, the operators ${\sf ADVEC}_s$ and {\sf ADVEC} should both return the mass flux divergence ${\sf DIVr}~{\tt v}$.

            \item {\bf Left null space}:
            
                 All constant fields should be in the null space of the adjoint operators ${\sf ADVEC}^*_s$ and ${\sf ADVEC}^*$.

            \item {\bf Meaningful adjoint}: 

                 ${\sf ADVEC}_s$'s adjoint is an accurate approximation of the adjoint of the continuous advection operator, like {\sf ADVEC}'s adjoint.
         \end{itemize}
         Because of the differences between ${\sf ADVEC}_s$ and ${\sf ADVEC}$ in the input and output spaces, certain interpolations and rotations are necessary which are not necessary in \cite{Win18WGW}, \cite{Win17WGK} and \cite{Hof12V}.}{Reviewer 1, Remark 5: Literature}

         The scalar advection operator ${\sf ADVEC}_s$ is defined as (cf. equation (\ref{eq: discrete DIV}))
         \begin{align}
            {\sf ADVEC}_s~{\tt f}_v :=
         \frac 12 \diag({\tt dVc})^{-1}
         \sum_{\eps=1/2,~3/2,~\ldots} \hspace{-0.5cm} & \left(
                                                      {\sf DIFFX}(\eps) \diag({\tt rho}_e(\eps)) \diag({\tt f}_e)\overline {\tt F}_e(\eps) \right. \nonumber \\
            &+
                                                    \left.  {\sf DIFFY}(\eps) \diag({\tt rho}_n(\eps))\diag({\tt f}_n)\overline {\tt F}_n(\eps) \right) \alpha(\eps).
         \label{eq: scalar ADVEC}
         \end{align}
         In this computation, $\overline{\tt F}_e$ and $\overline{\tt F}_n$ are computed using $\vec v$ and equations (\ref{eqn:interpolation_v}) and (\ref{eqn:matrixFen}).
         This advection operator ${\sf ADVEC}_s$ returns values at the pressure points of the grid. It satisfies the left null space property and the 'advection of a constant function is the divergence'-property
         \begin{eqnarray}
            \langle{\tt c1},{\sf ADVEC}_s~{\tt f}_v\rangle_c=0
            & \mbox{ and }  &
            {\sf ADVEC}_s~{\tt c1} =  {\sf DIVr}~{\tt v}.
         \end{eqnarray}
         The 'advection of a constant function is the divergence'-property is closely related to {\sf ADVEC}'s symmetry property in Table~\ref{Tab: properties}.

       In combination with the interpolation matrices {\sf C2E} and {\sf C2N}, the scalar advection operator can be used in models in which a scalar quantity is advected or convected, like temperature, or the concentration of a dissolved substance. The models discussed in this paper do not have such an advected/convected quantity.

       \bigskip
       The advection needed in the momentum equation of the isentropic Euler equations is not applied to a scalar field ${\tt f}_v$, but to a vector field ${\tt w}$. This is done by applying the scalar advection operator to the vector field's components separately. To make sure the operator is momentum preserving, this should not be done in the local grid orientation, but in the Cartesian orientation. The discretization process below is described in detail for 2D for simplicity, but all steps have also been worked out for 3D. 
         \begin{enumerate}
          \item cell-face interpolation and transformation to Cartesian orientation.
          
          The interpolation operators {\sf E2N} and {\sf N2E} of Section \ref{Sec: div at point i} are needed to calculate the complete vectors, represented by {\tt wx}, {\tt wy\_at\_e}, {\tt wx\_at\_n} and {\tt wy}.  Next, the vector field is transformed into the Cartesian grid. As an example, the vector field $\vec w_{e,i} = \vec r_x(\vec x_{e,i}) {\tt wx}_i + \vec r_y(\vec x_{e,i}) {\tt wy\_at\_e}_i$ is written as follows: 
          \[
           \vec w_{e,i} = \begin{pmatrix}
                           {\tt rxx\_at\_e}_i \\
                           {\tt rxy\_at\_e}_i \\
                          \end{pmatrix} {\tt wx}_i + 
                          \begin{pmatrix}
                           {\tt ryx\_at\_e}_i \\
                           {\tt ryy\_at\_e}_i \\
                          \end{pmatrix} {\tt wy\_at\_e}_i
                        = \begin{pmatrix}
                           {\tt rxx\_at\_e}_i {\tt wx}_i + {\tt ryx\_at\_e}_i {\tt wy\_at\_e}_i  \\
                           {\tt rxy\_at\_e}_i {\tt wx}_i + {\tt ryy\_at\_e}_i {\tt wy\_at\_e}_i  \\
                          \end{pmatrix}.
          \]
          This means that the component of $\vec w_{e,i}$ in the Cartesian direction $(1,0)$ equals
          \[
           {\tt w}_{e,(1,0),i} = (1,0) \cdot \vec w_{e,i} = {\tt rxx\_at\_e}_i {\tt wx}_i + {\tt ryx\_at\_e}_i {\tt wy\_at\_e}_i,
          \]
          such that
          \[
           {\tt w}_{e,(1,0)} = \diag({\tt rxx\_at\_e}) {\tt wx} + \diag({\tt ryx\_at\_e}) {\tt wy\_at\_e}.
          \]
          Similarly, the representation of $\vec w_{n,i}$ in the Cartesian basis vector $(1,0)$ can be computed:
          \[
           {\tt w}_{n,(1,0)} = \diag({\tt rxx\_at\_n}) {\tt wx\_at\_n} + \diag({\tt ryx\_at\_n}) {\tt wy}.
          \]
          The scalar field ${\tt w}_{v,(1,0)} = ({\tt w}_{e,(1,0)}^\top, {\tt w}_{n,(1,0)}^\top)^\top$ combines the two results. This scalar field has values in the complete velocity grid, like the scalar field ${\tt f}_v$ to which the advection operator ${\sf ADVEC}_s$ was applied earlier. It contains the horizontal component of the vector field ${\tt w}$.

          Using the same procedure, the vertical component is found:
         \[
           {\tt w}_{v,(0,1)} = \left(\begin{array}{l} {\tt w}_{e,(0,1)} \\ {\tt w}_{n,(0,1)}\end{array}\right) := \left( \begin{array}{c}
                   \diag({\tt rxy\_at\_e}) {\tt wx} + \diag({\tt ryy\_at\_e}) {\tt wy\_at\_e}\\
                   \diag({\tt rxy\_at\_n}) {\tt wx\_at\_n} + \diag({\tt ryy\_at\_n}) {\tt wy}\end{array}\right).
         \]

         \item ${\sf ADVEC}_c:$ application of ${\sf ADVEC}_s$ to each component.
         
         The advection operator ${\sf ADVEC}_s$ can be applied to ${\tt w}_{v,(1,0)}$ and ${\tt w}_{v,(0,1)}$, and the result will be the operator ${\sf ADVEC}_c$ (advection at pressure points):
         \[
            {\sf ADVEC}_c ~ {\tt w} := \left( \begin{array}{cc}
            {\sf ADVEC}_s ~ {\tt w}_{v,(1,0)} \\
            {\sf ADVEC}_s ~ {\tt w}_{v,(0,1)} \end{array}\right) :=
            \begin{pmatrix}
             {\tt a}_{c,(1,0)} \\
             {\tt a}_{c,(0,1)} 
            \end{pmatrix}
            .
         \]
       This operator approximates the continuous advection operator: $\vec a(\vec x_{c,i}) \hspace{-0.1cm} \approx ( ({\tt a}_{c,(1,0)})_i, ({\tt a}_{c,(0,1)})_i)$.
       This operator satisfies the left null space property given by
       \[
          \left\langle{\tt c1}, {\tt a}_{c,(1,0)} \right\rangle_c =
          \left\langle{\tt c1}, {\tt a}_{c,(0,1)} \right\rangle_c =0,
       \]
       and the 'advection of a constant function is the divergence'-property given by
       \[
          {\sf ADVEC}_c ({\tt c10},{\tt c01}) = \left( \begin{array}{cc} {\sf DIVr}~{\tt v} & {\sf 0} \\ {\sf 0} & {\sf DIVr}~{\tt v} \end{array}\right).
       \]

           \item ${\sf ADVEC}_a$: interpolation to the velocity grid and transformation to the curvilinear grid.
           
           The result of the advection operator ${\sf ADVEC}_c$ consists of standard Cartesian components at the pressure points of the grid, whereas the advection operator is needed for the components in local grid orientation at the velocity grid points. For this, we need to interpolate the components and transform the results. Interpolation is done by applying the adjoint operators ${\sf E2C}^*$ and ${\sf N2C}^*$, such that for each grid point $\vec x_{c,i}$, we compute
           \[ \vec a(\vec x_{e,i}) \approx
            \begin{pmatrix}
             ({\sf E2C}^* {\tt a}_{c,(1,0)})_i \\
             ({\sf E2C}^* {\tt a}_{c,(0,1)})_i
            \end{pmatrix}, \quad  \vec a(\vec x_{n,i}) \approx 
            \begin{pmatrix}
             ({\sf N2C}^* {\tt a}_{c,(1,0)})_i \\
             ({\sf N2C}^* {\tt a}_{c,(0,1)})_i
            \end{pmatrix}.
           \]
           Note that the interpolation from the cell centers to the cell faces is done with the adjoint interpolations ${\sf E2C}^*$ and ${\sf N2C}^*$ and not with the forward interpolations ${\sf C2E}$ and ${\sf C2N}$, as might have been expected. This is needed in order to satisfy the left null space property (\ref{eqn:fv_adveca}). In general, the adjoint of an accurate interpolation cannot be expected to be an accurate interpolation, but in the current case, the interpolations ${\sf E2C}$ and ${\sf N2C}$ are standard destaggering operations working on the computational grid, and their adjoints are accurate interpolations.
           
           Finally, the results are transformed to the basis spanned by the local grid orientation vectors: inner products with $r_x(\vec x_{e,i})$ and $r_y(\vec x_{n,i})$ are taken. Doing so defines the advection operator ${\sf ADVEC}_a$, which has almost all the properties listed in Table~\ref{Tab: properties}:
         \begin{eqnarray}
            {\sf ADVEC}_a = \left( \begin{array}{cc}
                   \diag({\tt rxx\_at\_e})~{\sf E2C}^* & \diag({\tt rxy\_at\_e})~{\sf E2C}^*\\
                   \diag({\tt ryx\_at\_n})~{\sf N2C}^* & \diag({\tt ryy\_at\_n})~{\sf N2C}^*\end{array}\right)
               ~{\sf ADVEC}_c.
         \end{eqnarray}
         
         In Section \ref{Sec: div at point i}, the interpolations {\sf E2C} and {\sf N2C} were required to be exact for the constant vector fields {\tt c10} and {\tt c01}.  The construction of the advection operator ${\sf ADVEC}_a$ also requires products of such fields to be exact, so (see equation (\ref{eqn:c100}) for the definition of {\tt c10} and {\tt c01})
       \begin{eqnarray*}
           {\sf E2C}~{\tt rxx\_at\_e} =  {\tt rxx\_at\_c} &,&
           {\sf E2C}~{\tt rxy\_at\_e} =  {\tt rxy\_at\_c}, \nonumber \\
           {\sf E2C}~\diag({\tt rxx\_at\_e}) {\tt rxx\_at\_e} &=& \diag({\tt rxx\_at\_c}) {\tt rxx\_at\_c}, \nonumber \\
           {\sf E2C}~\diag({\tt rxx\_at\_e}) {\tt rxy\_at\_e} &=& \diag({\tt rxx\_at\_c}) {\tt rxy\_at\_c}, \nonumber \\
           {\sf E2C}~\diag({\tt rxy\_at\_e}) {\tt rxy\_at\_e} &=& \diag({\tt rxy\_at\_c}) {\tt rxy\_at\_c},
       \end{eqnarray*}
         and similar for {\sf N2C}.

         This advection operator ${\sf ADVEC}_a$ also satisfies the left null space property. This can be seen by computing
         \begin{align*}
               \langle{\tt c10},{\sf ADVEC}_a~{\tt w}\rangle_v &= 
               \langle{\tt rxx\_at\_e}, \diag({\tt rxx\_at\_e})~{\sf E2C}^* {\sf ADVEC}_s~{\tt w}_{v,(1,0)}\rangle_e
                \\ &+
               \langle{\tt rxx\_at\_e}, \diag({\tt rxy\_at\_e})~{\sf E2C}^* {\sf ADVEC}_s~{\tt w}_{v,(0,1)}\rangle_e
                \\ &+
               \langle{\tt ryx\_at\_n}, \diag({\tt ryx\_at\_n})~{\sf N2C}^* {\sf ADVEC}_s~{\tt w}_{v,(1,0)} \rangle_n
                \\ &+
               \langle{\tt ryx\_at\_n}, \diag({\tt ryy\_at\_n})~{\sf N2C}^* {\sf ADVEC}_s~{\tt w}_{v,(0,1)}\rangle_n \\ &=
               \langle  \diag({\tt rxx\_at\_c}){\tt rxx\_at\_c}
                      + \diag({\tt ryx\_at\_c}){\tt ryx\_at\_c} , {\sf ADVEC}_s~{\tt w}_{v,(1,0)} \rangle_c
                \\ &+
               \langle    \diag({\tt rxy\_at\_c}){\tt rxx\_at\_c}
                       +  \diag({\tt ryy\_at\_c}){\tt ryx\_at\_c} , {\sf ADVEC}_s~{\tt w}_{v,(0,1)}\rangle_c.
         \end{align*}
         Due to orthonormality of the local grid orientation, we have
         \begin{align*}
          \begin{pmatrix}
           {\tt rxx\_at\_c} & {\tt rxy\_at\_c}  \\ 
           {\tt ryx\_at\_c} & {\tt ryy\_at\_c}
          \end{pmatrix}&
          \begin{pmatrix}
           {\tt rxx\_at\_c} & {\tt ryx\_at\_c} \\
           {\tt rxy\_at\_c} & {\tt ryy\_at\_c}
          \end{pmatrix}
        = 
        \begin{pmatrix}
         1 & 0 \\
         0 & 1 
        \end{pmatrix} \\
        &= 
          \begin{pmatrix}
           {\tt rxx\_at\_c} & {\tt ryx\_at\_c} \\
           {\tt rxy\_at\_c} & {\tt ryy\_at\_c}
          \end{pmatrix}
          \begin{pmatrix}
           {\tt rxx\_at\_c} & {\tt rxy\_at\_c}  \\ 
           {\tt ryx\_at\_c} & {\tt ryy\_at\_c}
          \end{pmatrix}
         ,
         \end{align*}
        such that
         \[
                {\tt rxx\_at\_c}_i^2 + {\tt ryx\_at\_c}_i^2 = 1, \quad
					  {\tt rxx\_at\_c}_i {\tt rxy\_at\_c}_i + {\tt ryx\_at\_c}_i {\tt ryy\_at\_c}_i = 0.
         \]
         Therefore, 
         \begin{equation}
          \langle{\tt c10},{\sf ADVEC}_a~{\tt w}\rangle_v = \langle  {\tt c1}, {\sf ADVEC}_s~{\tt w}_{v,(1,0)} \rangle_c  +
               \langle    0 , {\sf ADVEC}_s~{\tt w}_{v,(0,1)}\rangle_c = 0, \label{eqn:fv_adveca}
         \end{equation}
       and similar for $\langle {\tt c01}, {\sf ADVEC}_a~{\tt w}\rangle_v$.
                
       The advection operator ${\sf ADVEC}_a$ also has the following 'advection of a constant function is the divergence'-property:
         \[ {\sf ADVEC}_a~{\tt c} = \diag\left( \left( \begin{array}{cc} {\sf E2C}^* \\ {\sf N2C}^* \end{array} \right) {\sf DIVr}~{\tt v} \right) {\tt c}, \]
         for any discrete vector field {\tt c} that is the discrete representation of a constant vector field.

          \item {\sf ADVEC}: construction of the symmetry-preserving advection operator.
          
          The only property from Table~\ref{Tab: properties} that the advection operator ${\sf ADVEC}_a$ does not have, is the symmetry property.  We obtain the advection operator {\sf ADVEC} that has the symmetry property (as well as the left null space property) by combining the operator ${\sf ADVEC}_a$ and its adjoint \correction{\cite{Hof12V}}{Reviewer 1, Remark 5: Literature}:
         \[ {\sf ADVEC} := \frac 12 \left( {\sf ADVEC}_a - {\sf ADVEC}^*_a + \diag\left( \left( \begin{array}{cc} {\sf E2C}^* \\ {\sf N2C}^* \end{array} \right) {\sf DIVr}~{\tt v} \right) \right). \]
         The advection operator {\sf ADVEC} has precisely the symmetry property for the advection operator, provided that the interpolation operator ${\sf Interp}_{v \leftarrow c}$ used in the table is given by
         \[ {\sf Interp}_{v\leftarrow c} :=  \left( \begin{array}{cc} {\sf E2C}^* \\ {\sf N2C}^* \end{array} \right). \]

         \correction{Negative densities may be the result of this interpolation, in cases where the density rapidly approaches zero, for instance in cases of wetting 
         and drying, and this may affect the accuracy and stability of the numerical scheme.  In such cases, the order of interpolation and discretization 
         may be locally reduced (in a symmetry-preserving way). Alternatively, techniques like artificial porosity may be used to regularize the solution 
         \cite{Hof05V}. Such situations did not occur in the examples shown in this paper, and the techniques mentioned (local order reduction 
         and artificial porosity) are left outside the scope of this paper.}{Reviewer 1, Remark 7: positive rho?}
         \end{enumerate}
All required symmetry-preserving discrete operators are now defined, and we can investigate the performance of these operators for the three models. This will be done in the next section.

\section{Numerical results}\label{sec:results}
  In this section, the numerical results for the three different models are studied. The calculations for the linear-wave equations are done both with the Galerkin-type approach (Section \ref{sec:Galerkin}) and with the finite-difference technique (Section \ref{sec:FD}). This makes it possible to compare the results and computation times. The more advanced models for compressible waves and shallow waves are only investigated with the finite-difference technique.
  
\begin{figure}[ht!]
   \centering
            \subfigure[Uniform mesh]{\hspace*{-0.6cm}\includegraphics[scale = 0.4]{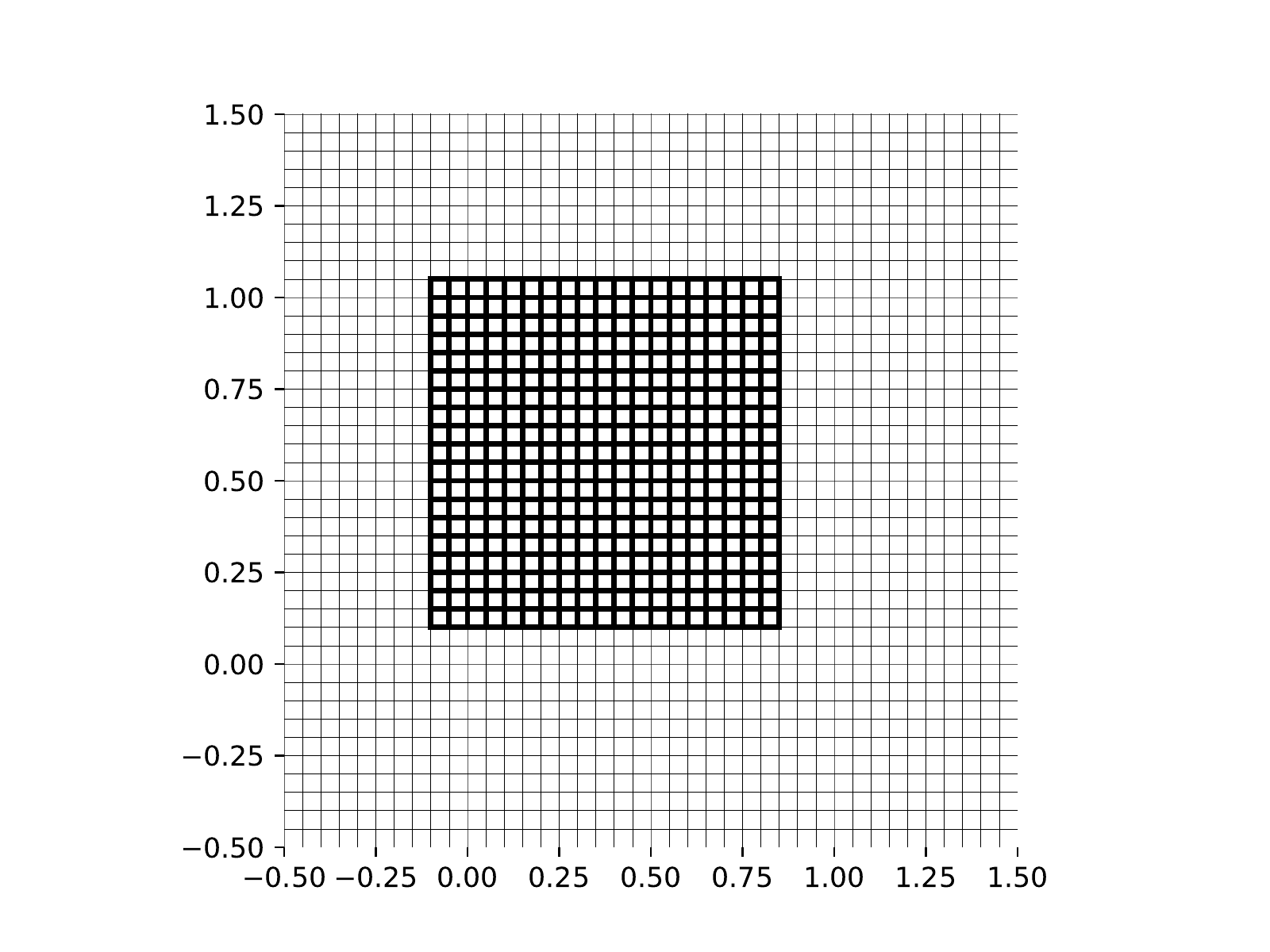}\hspace*{-1cm}}
            \subfigure[Curvilinear mesh]{\hspace*{-0.6cm}\includegraphics[scale = 0.4]{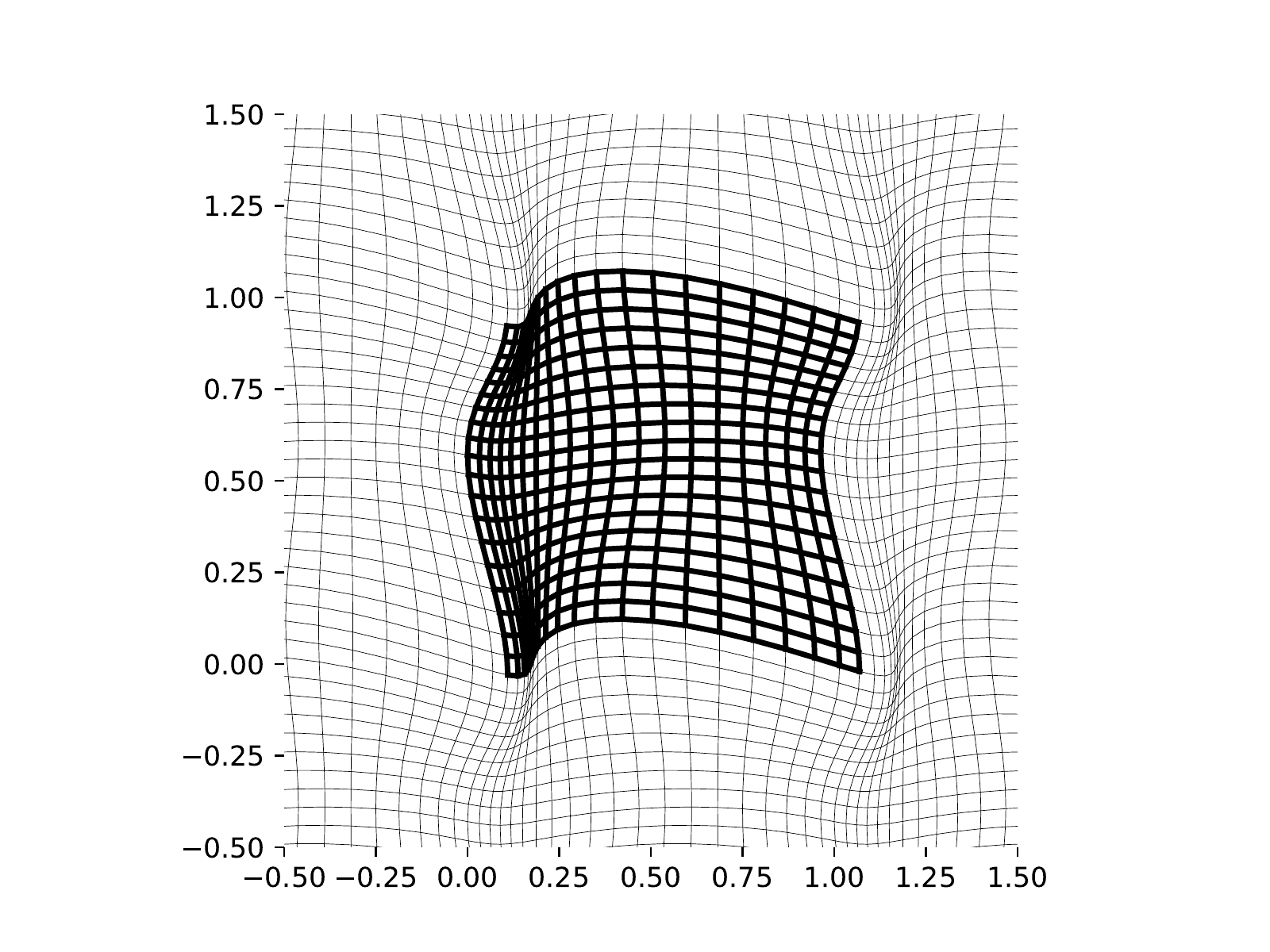}\hspace*{-1cm}}
            \subfigure[Detail of curvilinear mesh]{\hspace*{-0.6cm}\includegraphics[scale = 0.4]{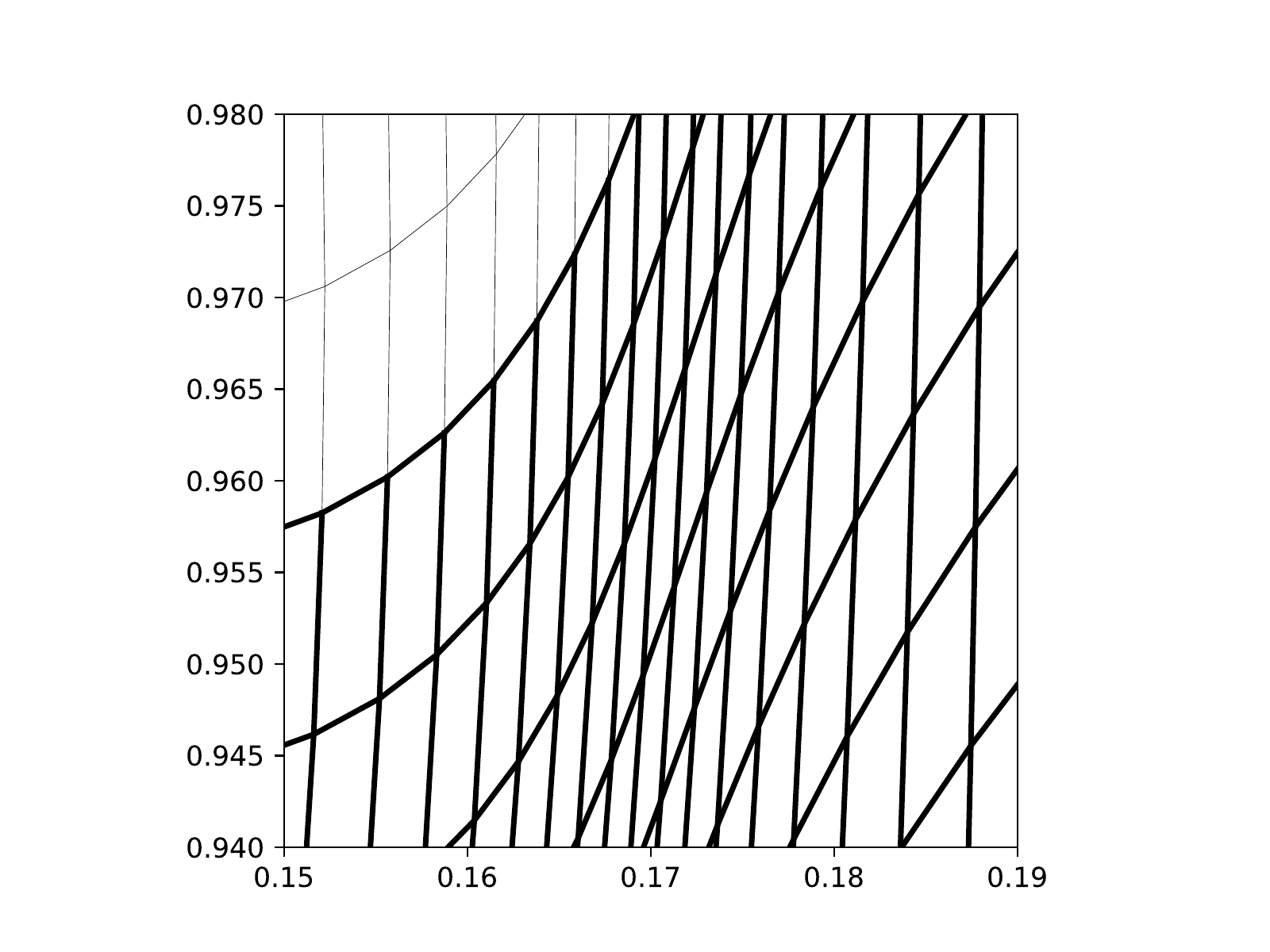}\hspace*{-1cm}}
            \caption{2D grids used in the experiments. An impression of the periodicity is given by extending the mesh on each boundary. 
                     \corr{The detail shows that the grid is non-orthogonal: in parts of the grid, the angle between grid lines are as small as 15$^{\mbox{\scriptsize o}}$}{Reviewer 2, Remarks 9 and 11: non-orthogonal grids}.}
         \label{Fig: grid}
         \end{figure}

  All calculations are done both on a uniform and on a periodic, non-orthogonal, curvilinear mesh. The periodicity of the curvilinear mesh gives us the opportunity to study symmetry preservation without focusing on boundary conditions (Section \ref{sec:properties}). An example of the uniform and curvilinear grids used for $20\times20$ grid points is given in Figure \ref{Fig: grid}.
  
   \subsection{Linear-wave equations}
      Choosing $p = \rho$ ($c=1$ in equation (\ref{eqn:linear-wave equations})), 
      an exact solution for the linear-wave equation is discussed in Section \ref{Sec: 4 steps}. Since the propagation speed $V_+=1$ only has one value, the solution is a traveling wave, which does not change its shape, but only its position, traveling in 'southeast' direction through the infinite periodic domain, given by
      \begin{eqnarray}
           p(x,y,t) = \exp\left(  -\frac{2}{9\pi^2}\sin^2(\sqrt 2\pi(x-y-\sqrt2 t))\right) &,&
       \vec v(x,y,t) = \frac 12 \sqrt 2 \left( \begin{array}r 1\\-1\end{array}\right) p(x,y,t).
      \end{eqnarray}

\begin{figure}[ht!]
         \centering
            \includegraphics[scale = 0.4]{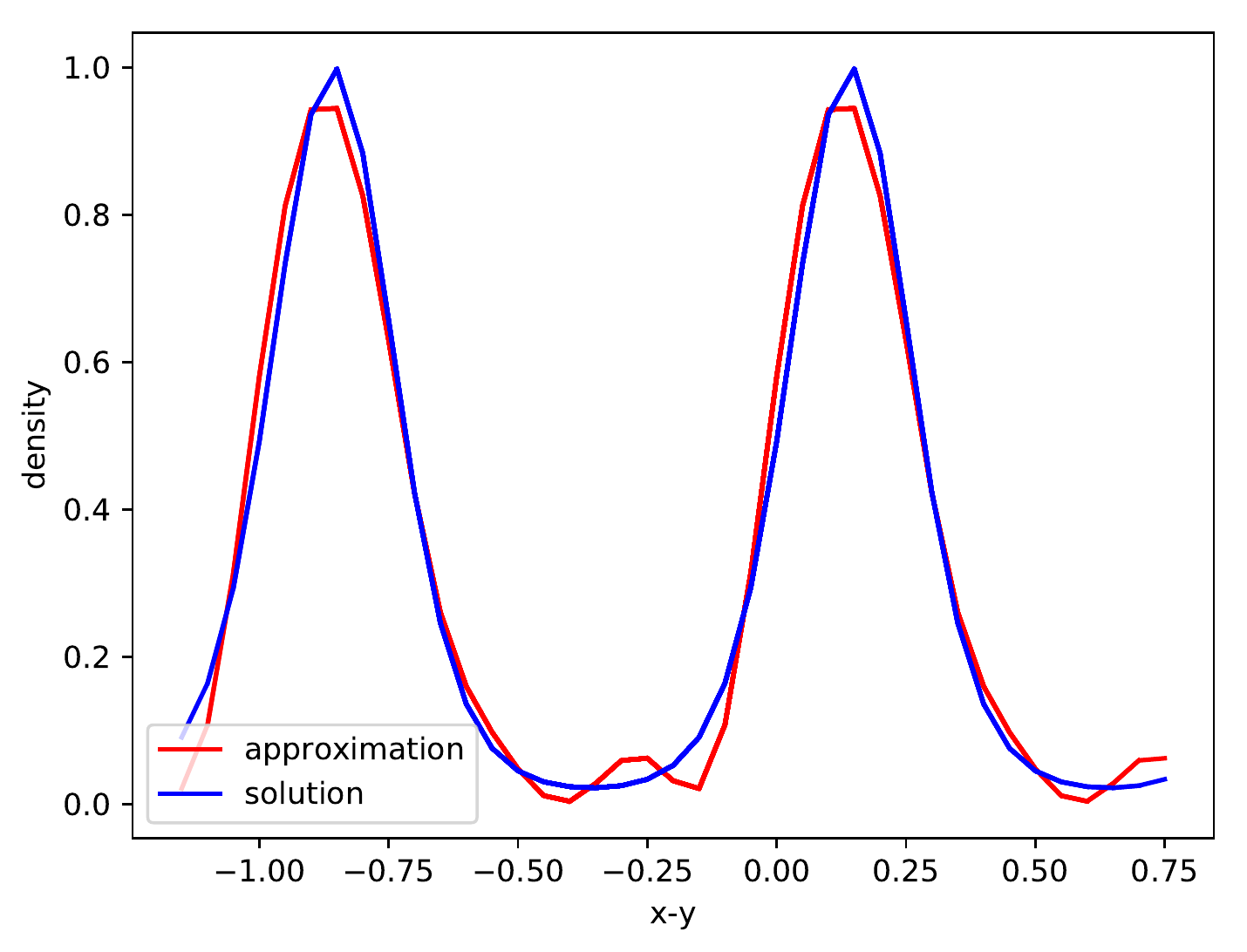} \label{Fig: linearFD res unif}
            \caption{Linear-wave equation: exact solution and numerical approximation of the density on a uniform mesh at $T=10$ (1D solution). The interpolation order is 4, $20\times20$ grid points.}
         \label{Fig: solutions linear wave eq}
         \end{figure}

         In fact, the solution on the uniform mesh is a 1D function of $x-y$, so the results at time $T=10$ can be represented in a 1D plot (Figure \ref{Fig: solutions linear wave eq}). The numerical approximation, shown in red in the figure, is close to the exact solution, shown in blue. The relative errors of the Galerkin-type approximations are shown in Table \ref{Tab: linear wave results} for different grid sizes and interpolation orders, and for the finite-difference approach they are shown in Table \ref{Tab: FD results}. It is seen that the solution becomes more accurate as more grid points are used and/or higher-order interpolation are used. The order of the approximation is equal to the interpolation order. Note that the results for our curvilinear grid are about the same as for the uniform grid. This is expected, since the exact solution, as well as the grid transformation $\vec X$ are infinitely smooth. 

         The effect of symmetry preservation is, that, assuming exact time integration, the total discrete mass {\tt M}, momentum $\vec{\tt M}$ and \correction{total energy}{Reviewer 2, Remark 7: energy equation} {\tt E} are constant up to machine precision, even in the most inaccurate solutions shown in Table \ref{Tab: linear wave results}. In Tables \ref{tab:conserve_linearGalerkin} and \ref{tab:conserve_FDGalerkin}, the mass, momentum and energy losses are given for a fourth-order discretization on a $20\times20$ grid.  \correction{Independently of the time-integration tolerance used, the mass is always conserved up to a very high precision, which depends on how accurately the matrix entries were calculated (using Richardson extrapolation of the composite trapezoidal rule).  For the finite-difference models, where the discrete operators require no numerical integration, the error for mass and momentum is really in the order of machine accuracy. Considering energy conservation, the error depends on the accuracy of the time integration.}{Reviewer 1, Remark 10: machine accuracy}

         \correction{
         For linear models, the implicit midpoint rule and its higher-order Gauss-Runge-Kutta generalizations will conserve energy (in fact all quadratic and linear invariants) up to machine precision.  For the nonlinear models, dedicated time-integration methods are needed. The current paper leaves time-integration methods outside the scope, because a standard time-integration method is used with such a small time step, that the time-integration errors are negligible.
         }{Reviewer 1, Remark 8: symplectic schemes}
                  
      Discrete mass, momentum and total energy are conserved up to machine precision, even in the most inaccurate solutions. \correction{The losses for different time-integration tolerances for the Galerkin-type approach are shown in Table \ref{tab:conserve_linearGalerkin}.  Similar results obtained in the finite-difference approach are shown in Table \ref{tab:conserve_FDGalerkin}.}{Reviewer 2, Remarks 2 and 3: conservation for all examples?}
      
      Because both the Galerkin-type approach and the finite-difference technique are available for the linear-wave equations, their relative computational performance can be compared. The Fortran code was compiled using gfortran version 7.3.0, and simulations were run on a desktop computer running Ubuntu 18.04 on an Intel i5-7400 CPU with 4 cores. Both methods use the same shared memory parallelization technique (OpenMP). We compare the results for interpolation (and derivative) order 8 and $40\times40$ grid points. 
      
      The Galerkin-type approach needs 118 seconds for the computation, thereby doing 21016 time steps and computing 27751 time derivatives. This means that this implementation computes 178 time steps per second, and 234 time-derivative evaluations per second. 
      
      On the other hand, the finite-difference approach needs 14 seconds for the total computation, in which it computes 26901 time steps and 31199 time derivatives. This implementation is able to compute 1922 time steps per second, and 2230 evaluations per second. We clearly see that the finite-difference approach is about ten times faster than the Galerkin-type approach. This also corresponds to the theory (Sections \ref{sec:computation_Galerkin} and \ref{Sec: div at point i}): for this 2D problem with 8th order interpolation, the Galerkin-type approach results in 480 nonzeros per matrix row, whereas the finite-difference approach only uses 48 nonzeros per grid point (a factor 10 difference).

         \begin{table}[ht]
            \caption{Relative errors (in 2-norm) for pressure in the linear-wave equations (Galerkin-type approach) at time $T=10$, for various interpolation orders, using reltol equal to {\tt 1e-11} in {\tt lsode}. For $160\times 160$ grid points on the uniform grid and 8th order, reltol={\tt 1e-13} is used.}
            \begin{center}
               \begin{tabular}{c|rrrlrlrlrlrlrlrlrlrlrlrlrlrl}
                  \multicolumn{9}{c}{\bf Uniform grid}\\
                                   &\multicolumn{2}{c}{order=2 }    &\multicolumn{2}{c}{order=4 }    &\multicolumn{2}{c}{order=6 }    &\multicolumn{2}{c}{order=8 }    \\
                                   & error            & order       & error            & order       & error            & order       & error            & order       \\ \hline
                  $10 \times 10$   & {\tt   1.12e+00} &             & {\tt   7.64e-01} &             & {\tt   1.68e-01} &             & {\tt   4.71e-01} &             \\
                  $20 \times 20$   & {\tt   1.43e+00} & {\tt  -0.35}& {\tt   4.10e-01} & {\tt   0.90}& {\tt   1.43e-01} & {\tt   0.23}& {\tt   4.49e-02} & {\tt   3.39}\\
                  $40 \times 40$   & {\tt   4.95e-01} & {\tt   1.53}& {\tt   7.85e-02} & {\tt   2.38}& {\tt   6.86e-03} & {\tt   4.38}& {\tt   7.82e-04} & {\tt   5.84}\\
                  $80 \times 80$   & {\tt   4.02e-01} & {\tt   0.30}& {\tt   6.34e-03} & {\tt   3.63}& {\tt   1.22e-04} & {\tt   5.82}& {\tt   3.55e-06} & {\tt   7.78}\\
                  $160 \times 160$ & {\tt   1.46e-01} & {\tt   1.46}& {\tt   4.01e-04} & {\tt   3.98}& {\tt   1.93e-06} & {\tt   5.98}& {\tt   1.43e-08} & {\tt   7.96}\\
               \end{tabular}\\[2ex]
               \begin{tabular}{c|rrrlrlrlrlrlrlrlrlrlrlrlrlrl}
                  \multicolumn{9}{c}{\bf Curvilinear grid}\\
                                   &\multicolumn{2}{c}{order=2 }    &\multicolumn{2}{c}{order=4 }    &\multicolumn{2}{c}{order=6 }    &\multicolumn{2}{c}{order=8 }    \\
                                   & error            & order       & error            & order       & error            & order       & error            & order       \\ \hline
                  $10 \times 10$   & {\tt   8.23e-01} &             & {\tt   1.06e+00} &             & {\tt   6.83e-01} &             & {\tt   3.57e+00} &             \\
                  $20 \times 20$   & {\tt   9.76e-01} & {\tt  -0.25}& {\tt   2.53e-01} & {\tt   2.07}& {\tt   4.29e-01} & {\tt   0.67}& {\tt   1.73e-01} & {\tt   4.36}\\
                  $40 \times 40$   & {\tt   8.67e-01} & {\tt   0.17}& {\tt   1.85e-01} & {\tt   0.45}& {\tt   4.61e-02} & {\tt   3.22}& {\tt   1.73e-02} & {\tt   3.33}\\
                  $80 \times 80$   & {\tt   5.22e-01} & {\tt   0.73}& {\tt   2.31e-02} & {\tt   3.00}& {\tt   1.49e-03} & {\tt   4.95}& {\tt   1.60e-04} & {\tt   6.75}\\
                  $160 \times 160$ & {\tt   2.03e-01} & {\tt   1.36}& {\tt   1.57e-03} & {\tt   3.88}& {\tt   2.51e-05} & {\tt   5.89}& {\tt   7.17e-06} & {\tt   7.80}\\
               \end{tabular}
            \end{center}
         \label{Tab: linear wave results}
         \end{table}

         \begin{table}[ht]
            \begin{center}
               \caption{Relative losses for linear-wave equations (Galerkin-type approach) at time $T=10$, for varying tolerances reltol in the time-integration method {\tt lsode}.
                        Fourth-order discretization on a $20\times 20$ grid.}\label{tab:conserve_linearGalerkin}
               \begin{tabular}{r|rrr|rrrrrrrrrr}
                      \multicolumn{1}{c}{} & \multicolumn{3}{c}{\bf Uniform grid} & \multicolumn{3}{c}{\bf Curvilinear grid}\\
                   {\bf reltol} &  {\bf mass loss} & {\bf mom. loss} & {\bf energy loss} &  {\bf mass loss} & {\bf mom. loss} & {\bf energy loss} \\ \hline
          \ignore{  {\tt 1e-06} &   {\tt 5.23e-14} &      {\tt 1.20e-13} &    {\tt 4.56e-05} &   {\tt 4.14e-15} &      {\tt 6.10e-15} &    {\tt 9.84e-06} \\}
                    {\tt 1e-07} &   {\tt 5.29e-14} &      {\tt 1.20e-13} &    {\tt 9.94e-07} &   {\tt 4.14e-15} &      {\tt 6.33e-15} &    {\tt 1.60e-06} \\
                    {\tt 1e-08} &   {\tt 5.27e-14} &      {\tt 1.20e-13} &    {\tt 6.97e-08} &   {\tt 2.83e-15} &      {\tt 4.78e-15} &    {\tt 1.56e-08} \\
                    {\tt 1e-09} &   {\tt 5.16e-14} &      {\tt 1.21e-13} &    {\tt 6.39e-08} &   {\tt 3.92e-15} &      {\tt 6.57e-15} &    {\tt 3.25e-09} \\
                    {\tt 1e-10} &   {\tt 5.29e-14} &      {\tt 1.20e-13} &    {\tt 1.15e-09} &   {\tt 3.92e-15} &      {\tt 5.58e-15} &    {\tt 1.50e-11} \\
                    {\tt 1e-11} &   {\tt 5.36e-14} &      {\tt 1.17e-13} &    {\tt 5.54e-11} &   {\tt 2.61e-15} &      {\tt 7.73e-15} &    {\tt 1.87e-11} \\
                    {\tt 1e-12} &   {\tt 5.23e-14} &      {\tt 1.18e-13} &    {\tt 4.74e-11} &   {\tt 1.96e-15} &      {\tt 4.43e-15} &    {\tt 5.84e-12} \\
               \end{tabular}
             \end{center}
         \end{table}
         
         \begin{table}[ht]
            \caption{Relative errors (in 2-norm) for pressure in the linear-wave equations (finite-difference approach) at time $T=10$, for various (same) orders for derivatives and interpolations, using reltol equal to {\tt 1e-11} in {\tt lsode}. For $160\times 160$ grid points and 8th order, reltol={\tt 1e-13} is used.}
            \begin{center}
               \begin{tabular}{c|rlrlrlrr}
                  \multicolumn{9}{c}{\bf Uniform grid}\\
                                   &\multicolumn{2}{c}{order=2 }    &\multicolumn{2}{c}{order=4 }    &\multicolumn{2}{c}{order=6 }    &\multicolumn{2}{c}{order=8 }    \\
                                   & error            & order       & error            & order       & error            & order       & error            & order       \\ \hline
                  $10 \times 10$   & {\tt   9.29e-01} &             & {\tt   4.52e-01} &             & {\tt   1.33e-01} &             & {\tt   1.62e-01} &             \\
                  $20 \times 20$   & {\tt   6.22e-01} & {\tt   0.58}& {\tt   8.90e-02} & {\tt   2.34}& {\tt   1.88e-02} & {\tt   2.82}& {\tt   5.54e-03} & {\tt   4.87}\\
                  $40 \times 40$   & {\tt   2.58e-01} & {\tt   1.27}& {\tt   7.68e-03} & {\tt   3.53}& {\tt   4.14e-04} & {\tt   5.51}& {\tt   3.59e-05} & {\tt   7.27}\\
                  $80 \times 80$   & {\tt   7.67e-02} & {\tt   1.75}& {\tt   4.91e-04} & {\tt   3.97}& {\tt   6.77e-06} & {\tt   5.94}& {\tt   1.49e-07} & {\tt   7.92}\\
                  $160 \times 160$ & {\tt   1.96e-02} & {\tt   1.97}& {\tt   3.08e-05} & {\tt   4.00}& {\tt   9.79e-08} & {\tt   6.11}& {\tt   6.00e-10} & {\tt   7.96}\\
               \end{tabular}\\[2ex]
               \begin{tabular}{c|rlrlrlrlrlrlrlrlrlrlrlrlrlrr}
                  \multicolumn{9}{c}{\bf Curvilinear grid}\\
                                   &\multicolumn{2}{c}{order=2 }    &\multicolumn{2}{c}{order=4 }    &\multicolumn{2}{c}{order=6 }    &\multicolumn{2}{c}{order=8 }    \\
                                   & error            & order       & error            & order       & error            & order       & error            & order       \\ \hline
                  $10 \times 10$   & {\tt   1.14e+00} &             & {\tt   5.08e-01} &             & {\tt   2.92e-01} &             & {\tt   4.76e-01} &             \\
                  $20 \times 20$   & {\tt   4.93e-01} & {\tt   1.21}& {\tt   1.65e-01} & {\tt   1.62}& {\tt   6.65e-02} & {\tt   2.13}& {\tt   2.29e-02} & {\tt   4.37}\\
                  $40 \times 40$   & {\tt   3.31e-01} & {\tt   0.57}& {\tt   2.13e-02} & {\tt   2.95}& {\tt   2.58e-03} & {\tt   4.69}& {\tt   5.36e-04} & {\tt   5.42}\\
                  $80 \times 80$   & {\tt   1.13e-01} & {\tt   1.55}& {\tt   1.43e-03} & {\tt   3.90}& {\tt   4.47e-05} & {\tt   5.85}& {\tt   2.60e-06} & {\tt   7.69}\\
                  $160 \times 160$ & {\tt   3.03e-02} & {\tt   1.91}& {\tt   9.01e-05} & {\tt   3.99}& {\tt   7.17e-07} & {\tt   5.96}& {\tt   1.05e-08} & {\tt   7.95}\\
               \end{tabular}
            \end{center}
         \label{Tab: FD results}
         \end{table}

         \begin{table}[ht]
            \begin{center}
               \caption{Relative losses for linear-wave equations (finite-difference approach) at time $T=10$, for varying tolerances reltol in the time-integration method {\tt lsode}.
                        Fourth-order discretization on a $20\times 20$ grid.}\label{tab:conserve_FDGalerkin}
               \begin{tabular}{r|rrr|rrrrrrrrrr}
                      \multicolumn{1}{c}{} & \multicolumn{3}{c}{\bf Uniform grid} & \multicolumn{3}{c}{\bf Curvilinear grid}\\
                    {\bf reltol} &  {\bf mass loss} & {\bf mom. loss} & {\bf energy loss}&  {\bf mass loss} & {\bf mom. loss} & {\bf energy loss} \\ \hline
       \ignore{     {\tt 1e-05} &   {\tt 2.18e-16} &      {\tt 1.66e-15} &    {\tt 5.22e-05} &   {\tt 4.35e-16} &      {\tt 1.23e-14} &    {\tt 9.97e-05} \\
                    {\tt 1e-06} &   {\tt 8.71e-16} &      {\tt 1.57e-15} &    {\tt 3.53e-05} &   {\tt 2.18e-16} &      {\tt 1.41e-14} &    {\tt 6.46e-07} \\}
                    {\tt 1e-07} &   {\tt 2.18e-16} &      {\tt 2.39e-15} &    {\tt 1.94e-06} &   {\tt 2.18e-16} &      {\tt 1.36e-14} &    {\tt 3.29e-08} \\
                    {\tt 1e-08} &   {\tt 1.09e-15} &      {\tt 1.57e-15} &    {\tt 1.22e-06} &   {\tt 2.18e-16} &      {\tt 1.43e-14} &    {\tt 4.61e-09} \\
                    {\tt 1e-09} &   {\tt 4.35e-16} &      {\tt 1.12e-15} &    {\tt 1.78e-09} &   {\tt 4.35e-16} &      {\tt 1.15e-14} &    {\tt 4.33e-10} \\
                    {\tt 1e-10} &   {\tt 1.09e-15} &      {\tt 1.66e-15} &    {\tt 2.98e-09} &   {\tt 1.52e-15} &      {\tt 1.10e-14} &    {\tt 7.85e-11} \\
                    {\tt 1e-11} &   {\tt 0.00e+00} &      {\tt 9.86e-16} &    {\tt 1.35e-09} &   {\tt 6.53e-16} &      {\tt 1.52e-14} &    {\tt 1.28e-11} \\
                    {\tt 1e-12} &   {\tt 2.18e-16} &      {\tt 1.66e-15} &    {\tt 3.16e-11} &   {\tt 1.31e-15} &      {\tt 1.27e-14} &    {\tt 8.80e-13} \\
               \end{tabular}
            \end{center}
         \end{table}

   \subsection{Compressible-wave equations}
\begin{figure}[ht!]
      \centering
         \subfigure[Initial solution]{\includegraphics[scale = 0.4]{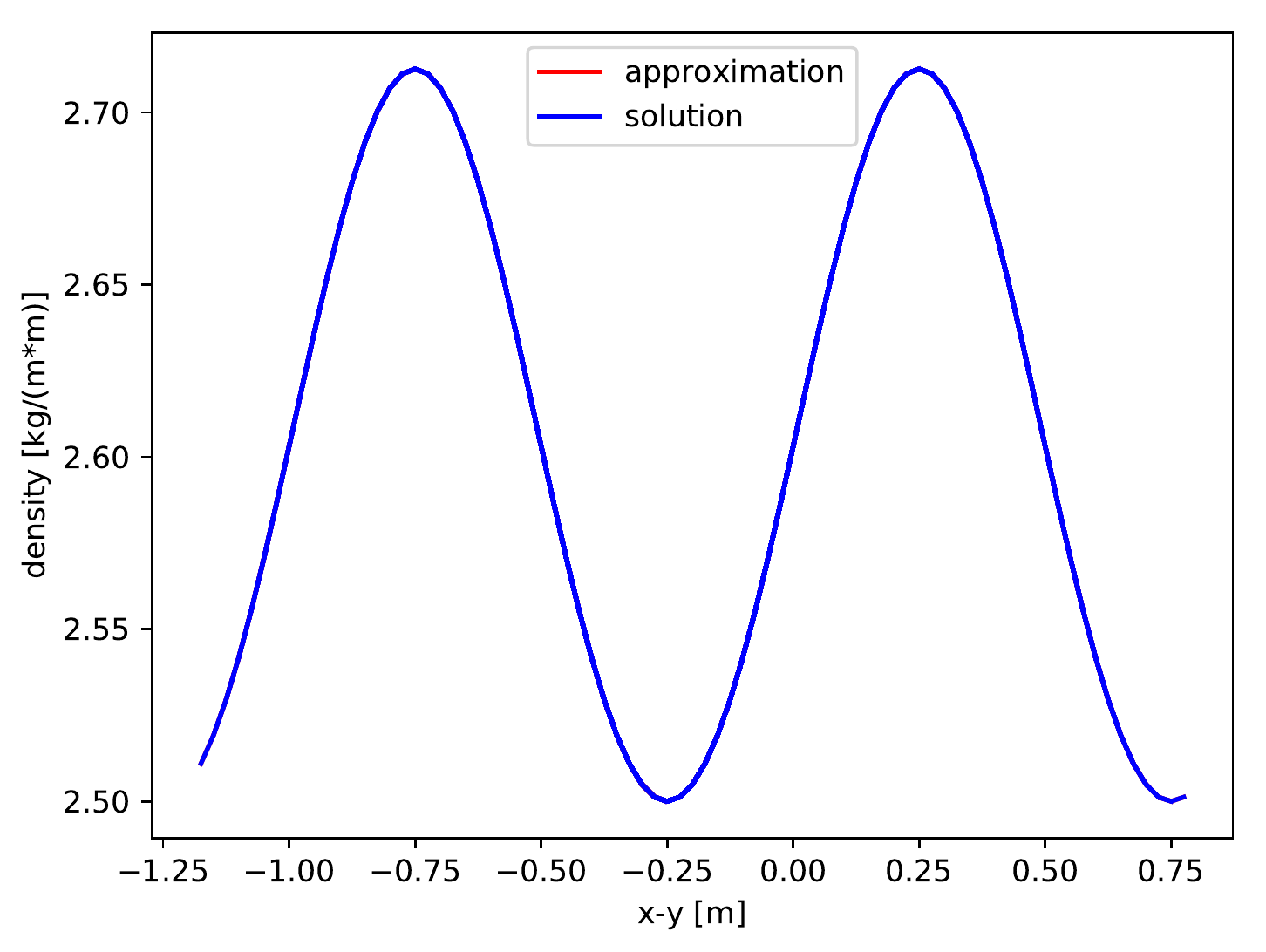}}
         \subfigure[Final-time approximation]{\includegraphics[scale = 0.4]{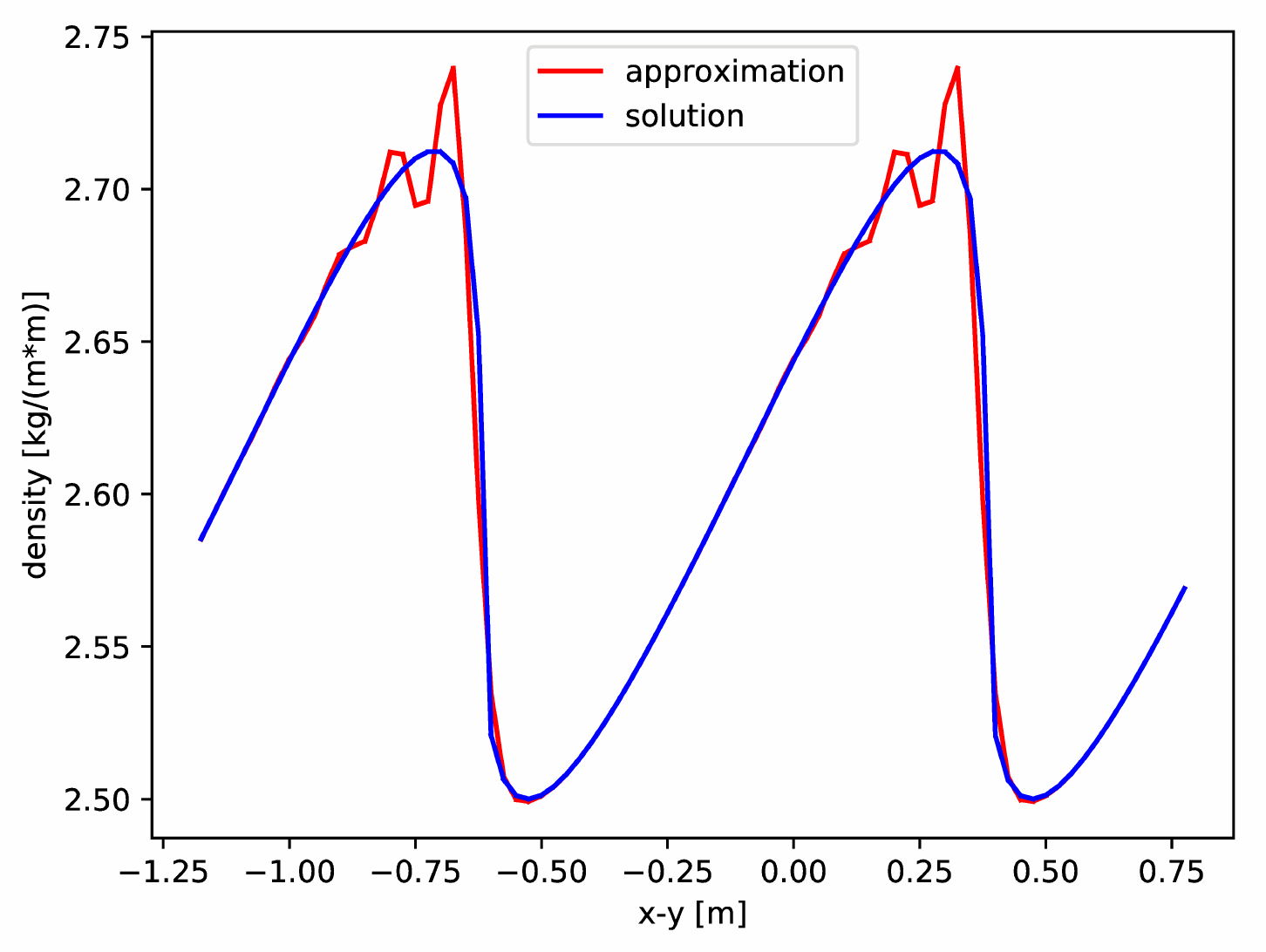}}
         \caption{Exact solution and numerical approximation for the compressible-wave equations on a uniform, $40\times40$ mesh (1D solution), with interpolation order 4. The final-time approximation belongs to the moment when the solution becomes discontinuous.}
      \label{Fig: solutions compressible wave eq}
      \end{figure}

      Unlike the solutions in the previous numerical examples, the solution of the compressible-wave equations has variable propagation speed. This means that the crests of the wave travel faster than the troughs.  This effect causes the wave to steepen over time, and finally to become discontinuous. In Figure \ref{Fig: solutions compressible wave eq}, the initial solution and the numerical approximation at the time $t_N$ (when the solution becomes discontinuous) are shown.  Obviously, the solution becomes more difficult to approximate for a numerical method as the moment $t_N$ approaches. However, we see that away from the shock, the approximations are still accurate.       Because the solution becomes increasingly difficult to approximate as time $t=t_N$ approaches, results are shown for $t=1.09=t_N/2$ and for $t=2.18=t_N$ in Table \ref{Tab: State results final}. Halfway the simulation, the convergence order approaches the theoretical order; for the final time, only order 1 can be achieved. 
      
      \begin{table}[p]
         \caption{Relative errors (in 2-norm) for the densities in the compressible-wave equations, for various (same) orders for derivatives and interpolations, using reltol equal to {\tt 1e-11} in {\tt lsode}.}
         \begin{center}
            \begin{tabular}{c|rlrlrlrlrlrlrlrlrlrlrlrlrlrl}
               \multicolumn{9}{c}{\bf Uniform grid, $t=t_N/2$}\\
                                &\multicolumn{2}{c}{order=2 }    &\multicolumn{2}{c}{order=4 }    &\multicolumn{2}{c}{order=6 }    &\multicolumn{2}{c}{order=8 }    \\
                                & error            & order       & error            & order       & error            & order       & error            & order       \\ \hline
               $10 \times 10$   & {\tt   3.03e-01} &             & {\tt   6.65e-02} &             & {\tt   3.52e-02} &             & {\tt   2.48e-02} &             \\
               $20 \times 20$   & {\tt   1.04e-01} & {\tt   1.54}& {\tt   1.24e-02} & {\tt   2.42}& {\tt   5.29e-03} & {\tt   2.74}& {\tt   3.49e-03} & {\tt   2.83}\\
               $40 \times 40$   & {\tt   3.05e-02} & {\tt   1.77}& {\tt   1.32e-03} & {\tt   3.24}& {\tt   3.13e-04} & {\tt   4.08}& {\tt   1.53e-04} & {\tt   4.51}\\
               $80 \times 80$   & {\tt   7.98e-03} & {\tt   1.93}& {\tt   9.50e-05} & {\tt   3.79}& {\tt   8.55e-06} & {\tt   5.19}& {\tt   2.29e-06} & {\tt   6.06}\\
               $160 \times 160$ & {\tt   2.01e-03} & {\tt   1.99}& {\tt   6.10e-06} & {\tt   3.96}& {\tt   1.59e-07} & {\tt   5.75}& {\tt   1.56e-08} & {\tt   7.20}\\
            \end{tabular}\\[2ex]
            \begin{tabular}{c|rlrlrlrlrlrlrlrlrlrlrlrlrlrl}
               \multicolumn{9}{c}{\bf Curvilinear grid, $t=t_N/2$}\\
                                &\multicolumn{2}{c}{order=2 }    &\multicolumn{2}{c}{order=4 }    &\multicolumn{2}{c}{order=6 }    &\multicolumn{2}{c}{order=8 }    \\
                                & error            & order       & error            & order       & error            & order       & error            & order       \\ \hline
               $10 \times 10$   & {\tt   4.35e-01} &             & {\tt   1.38e-01} &             & {\tt   9.38e-02} &             & {\tt   7.94e-02} &             \\
               $20 \times 20$   & {\tt   1.55e-01} & {\tt   1.48}& {\tt   2.95e-02} & {\tt   2.23}& {\tt   1.67e-02} & {\tt   2.49}& {\tt   1.26e-02} & {\tt   2.65}\\
               $40 \times 40$   & {\tt   4.75e-02} & {\tt   1.71}& {\tt   4.67e-03} & {\tt   2.66}& {\tt   1.90e-03} & {\tt   3.14}& {\tt   1.22e-03} & {\tt   3.37}\\
               $80 \times 80$   & {\tt   1.26e-02} & {\tt   1.91}& {\tt   4.14e-04} & {\tt   3.49}& {\tt   8.85e-05} & {\tt   4.43}& {\tt   4.10e-05} & {\tt   4.89}\\
               $160 \times 160$ & {\tt   3.15e-03} & {\tt   2.00}& {\tt   2.75e-05} & {\tt   3.91}& {\tt   2.12e-06} & {\tt   5.39}& {\tt   5.13e-07} & {\tt   6.32}\\
            \end{tabular}\\[2ex]
               \begin{tabular}{c|rlrlrlrlrlrlrlrlrlrlrlrlrlrl}
                  \multicolumn{9}{c}{\bf Uniform grid, $t=t_N$}\\
                                   &\multicolumn{2}{c}{order=2 }    &\multicolumn{2}{c}{order=4 }    &\multicolumn{2}{c}{order=6 }    &\multicolumn{2}{c}{order=8 }    \\
                                   & error            & order       & error            & order       & error            & order       & error            & order       \\ \hline
                  $10 \times 10$   & {\tt   5.42e-01} &             & {\tt   2.62e-01} &             & {\tt   2.20e-01} &             & {\tt   1.96e-01} &             \\
                  $20 \times 20$   & {\tt   3.25e-01} & {\tt   0.74}& {\tt   1.44e-01} & {\tt   0.86}& {\tt   1.16e-01} & {\tt   0.92}& {\tt   1.07e-01} & {\tt   0.87}\\
                  $40 \times 40$   & {\tt   1.82e-01} & {\tt   0.84}& {\tt   7.63e-02} & {\tt   0.92}& {\tt   5.92e-02} & {\tt   0.97}& {\tt   5.24e-02} & {\tt   1.03}\\
                  $80 \times 80$   & {\tt   1.04e-01} & {\tt   0.82}& {\tt   4.06e-02} & {\tt   0.91}& {\tt   3.14e-02} & {\tt   0.91}& {\tt   2.81e-02} & {\tt   0.90}\\
                  $160 \times 160$ & {\tt   6.09e-02} & {\tt   0.77}& {\tt   2.10e-02} & {\tt   0.95}& {\tt   1.49e-02} & {\tt   1.08}& {\tt   1.28e-02} & {\tt   1.14}\\
               \end{tabular}\\[2ex]
               \begin{tabular}{c|rlrlrlrlrlrlrlrlrlrlrlrlrlrl}
                  \multicolumn{9}{c}{\bf Curvilinear grid, $t=t_N$}\\
                                   &\multicolumn{2}{c}{order=2 }    &\multicolumn{2}{c}{order=4 }    &\multicolumn{2}{c}{order=6 }    &\multicolumn{2}{c}{order=8 }    \\
                                   & error            & order       & error            & order       & error            & order       & error            & order       \\ \hline
                  $10 \times 10$   & {\tt   7.41e-01} &             & {\tt   3.19e-01} &             & {\tt   2.24e-01} &             & {\tt   2.15e-01} &             \\
                  $20 \times 20$   & {\tt   4.13e-01} & {\tt   0.84}& {\tt   1.83e-01} & {\tt   0.80}& {\tt   1.43e-01} & {\tt   0.64}& {\tt   1.30e-01} & {\tt   0.73}\\
                  $40 \times 40$   & {\tt   2.33e-01} & {\tt   0.83}& {\tt   1.02e-01} & {\tt   0.84}& {\tt   8.10e-02} & {\tt   0.82}& {\tt   7.29e-02} & {\tt   0.83}\\
                  $80 \times 80$   & {\tt   1.39e-01} & {\tt   0.75}& {\tt   5.98e-02} & {\tt   0.77}& {\tt   4.63e-02} & {\tt   0.81}& {\tt   4.10e-02} & {\tt   0.83}\\
                  $160 \times 160$ & {\tt   8.34e-02} & {\tt   0.74}& {\tt   3.36e-02} & {\tt   0.83}& {\tt   2.51e-02} & {\tt   0.88}& {\tt   2.20e-02} & {\tt   0.90}\\
               \end{tabular}
         \end{center}
         \label{Tab: State results final}
      \end{table}

      \begin{table}[ht]
           \caption{\corr{Relative losses for the compressible-wave equations, for varying tolerances reltol in the time-integration method {\tt lsode}, at the time when the solution becomes discontinuous.
                    Fourth-order discretization on a $20\times 20$ grid.}{Reviewer 2, Remarks 2 and 3: conservation for all examples} }\label{tab:conservation_compressible}
       \begin{center}
           \begin{tabular}{r|rrr|rrrrrrrrrr}
                  \multicolumn{1}{c}{} & \multicolumn{3}{c}{\bf Uniform grid} & \multicolumn{3}{c}{\bf Curvilinear grid}\\
                 {\bf reltol} &  {\bf mass loss} & {\bf mom. loss} & {\bf energy loss} &  {\bf mass loss} & {\bf mom. loss} & {\bf energy loss} \\ \hline
                  {\tt 1e-05} &   {\tt 6.82e-16} &      {\tt 1.85e-14} &    {\tt 6.53e-09} &   {\tt 1.02e-15} &      {\tt 1.67e-13} &    {\tt 6.41e-10} \\
                  {\tt 1e-06} &   {\tt 3.41e-16} &      {\tt 1.94e-14} &    {\tt 7.99e-10} &   {\tt 1.19e-15} &      {\tt 1.69e-13} &    {\tt 1.57e-10} \\
                  {\tt 1e-07} &   {\tt 1.02e-15} &      {\tt 2.14e-14} &    {\tt 4.18e-11} &   {\tt 8.52e-16} &      {\tt 1.64e-13} &    {\tt 6.64e-12} \\
                  {\tt 1e-08} &   {\tt 3.41e-16} &      {\tt 2.39e-14} &    {\tt 1.68e-13} &   {\tt 0.00e+00} &      {\tt 1.68e-13} &    {\tt 6.49e-13} \\
                  {\tt 1e-09} &   {\tt 6.82e-16} &      {\tt 1.75e-14} &    {\tt 3.30e-13} &   {\tt 8.52e-16} &      {\tt 1.69e-13} &    {\tt 5.20e-14} \\
                  {\tt 1e-10} &   {\tt 1.70e-16} &      {\tt 1.81e-14} &    {\tt 3.33e-14} &   {\tt 5.11e-16} &      {\tt 1.66e-13} &    {\tt 1.09e-15} \\
                  {\tt 1e-11} &   {\tt 3.41e-16} &      {\tt 1.97e-14} &    {\tt 2.90e-15} &   {\tt 2.56e-15} &      {\tt 1.66e-13} &    {\tt 1.45e-15} \\
                  {\tt 1e-12} &   {\tt 1.02e-15} &      {\tt 2.10e-14} &    {\tt 1.81e-16} &   {\tt 1.88e-15} &      {\tt 1.66e-13} &    {\tt 1.81e-16} \\
           \end{tabular}
       \end{center}
      \end{table}
      
      \correction{Table \ref{tab:conservation_compressible} shows that even the most inaccurate solutions at the final time conserve discrete mass, momentum and energy up to machine precision.}{Reviewer 2, Remarks 2 and 3: conservation for all examples} However, the physically-relevant solution should have lower \correction{total energy}{Reviewer 2, Remark 7: energy equation}  after the shock has appeared. Therefore, numerical viscosity can be added, for instance in the form of flux limiters, but such techniques are beyond the scope of this paper.

   \subsection{Shallow-water equations}
      \begin{figure}[ht!]
         \centering
         \includegraphics[scale = 0.4]{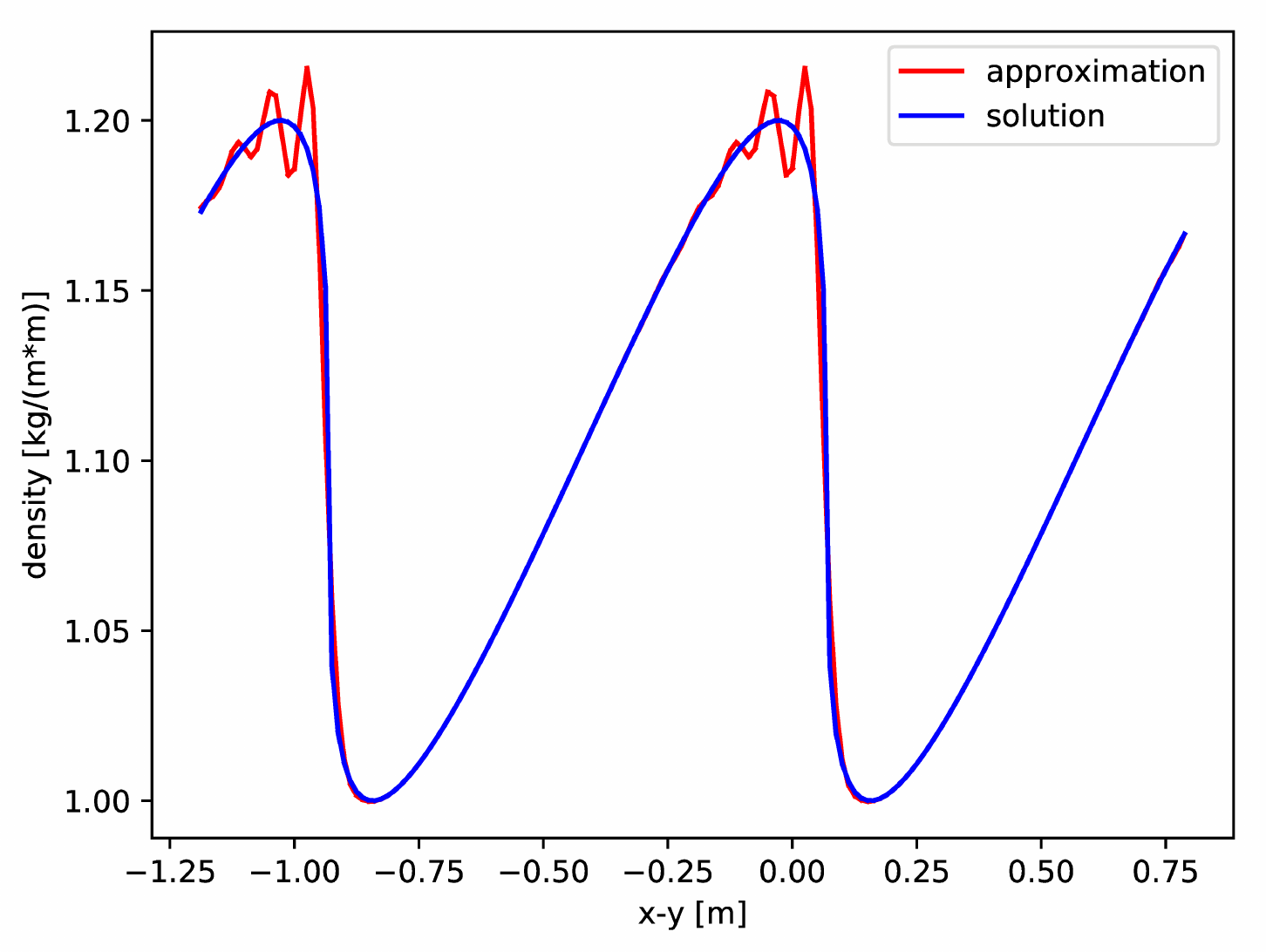}
         \caption{Shallow-water equations: exact solution and numerical approximation at $t=t_N$ on a uniform mesh (1D solution). The interpolation order is 4, and the meshes contain $80\times80$ grid points.}
      \label{Fig: solutions shallow}
      \end{figure}
      Like the compressible-wave equations, the variable propagation speed leads to a steepening wave in the exact solution. The results at time $t=t_N$ (the moment that the solution becomes discontinuous) are visualized in Figure \ref{Fig: solutions shallow}. The accuracies at time $t=t_N/2$ are shown in Table \ref{Tab: Shallow results final}: the accuracy orders are approaching the theoretical order. At the final time, only an order 1 can be achieved. \correction{Mass, momentum and energy are again conserved up to machine precision, as shown in Table \ref{tab:conservation_shallow}.}{Reviewer 2, Remarks 2 and 3: conservation for all examples}
      
      \begin{table}[ht]
         \caption{Relative errors (in 2-norm) for the densities in the shallow-water equations at time $t=t_N/2$, for various (same) orders for derivatives and interpolations, using reltol equal to {\tt 1e-11} in {\tt lsode}.}
         \begin{center}
               \begin{tabular}{c|rlrlrlrlrlrlrlrlrlrlrlrlrlrl}
                  \multicolumn{9}{c}{\bf Uniform grid}\\
                                   &\multicolumn{2}{c}{order=2 }    &\multicolumn{2}{c}{order=4 }    &\multicolumn{2}{c}{order=6 }    &\multicolumn{2}{c}{order=8 }    \\
                                   & error            & order       & error            & order       & error            & order       & error            & order       \\ \hline
                  $10 \times 10$   & {\tt   7.07e-02} &             & {\tt   2.23e-02} &             & {\tt   1.33e-02} &             & {\tt   1.10e-02} &             \\
                  $20 \times 20$   & {\tt   2.62e-02} & {\tt   1.43}& {\tt   4.60e-03} & {\tt   2.28}& {\tt   2.31e-03} & {\tt   2.53}& {\tt   1.65e-03} & {\tt   2.74}\\
                  $40 \times 40$   & {\tt   8.16e-03} & {\tt   1.68}& {\tt   6.00e-04} & {\tt   2.94}& {\tt   1.89e-04} & {\tt   3.61}& {\tt   1.05e-04} & {\tt   3.97}\\
                  $80 \times 80$   & {\tt   2.21e-03} & {\tt   1.89}& {\tt   4.81e-05} & {\tt   3.64}& {\tt   6.40e-06} & {\tt   4.89}& {\tt   2.03e-06} & {\tt   5.69}\\
                  $160 \times 160$ & {\tt   5.61e-04} & {\tt   1.98}& {\tt   3.15e-06} & {\tt   3.93}& {\tt   1.24e-07} & {\tt   5.69}& {\tt   1.46e-08} & {\tt   7.12}\\
               \end{tabular}\\[2ex]
               \begin{tabular}{c|rlrlrlrlrlrlrlrlrlrlrlrlrlrl}
                  \multicolumn{9}{c}{\bf Curvilinear grid}\\
                                   &\multicolumn{2}{c}{order=2 }    &\multicolumn{2}{c}{order=4 }    &\multicolumn{2}{c}{order=6 }    &\multicolumn{2}{c}{order=8 }    \\
                                   & error            & order       & error            & order       & error            & order       & error            & order       \\ \hline
                  $10 \times 10$   & {\tt   8.48e-02} &             & {\tt   4.35e-02} &             & {\tt   3.51e-02} &             & {\tt   3.49e-02} &             \\
                  $20 \times 20$   & {\tt   3.34e-02} & {\tt   1.34}& {\tt   1.17e-02} & {\tt   1.90}& {\tt   7.73e-03} & {\tt   2.18}& {\tt   6.16e-03} & {\tt   2.50}\\
                  $40 \times 40$   & {\tt   1.25e-02} & {\tt   1.42}& {\tt   2.19e-03} & {\tt   2.42}& {\tt   1.06e-03} & {\tt   2.87}& {\tt   7.32e-04} & {\tt   3.07}\\
                  $80 \times 80$   & {\tt   3.71e-03} & {\tt   1.75}& {\tt   2.34e-04} & {\tt   3.22}& {\tt   6.26e-05} & {\tt   4.08}& {\tt   3.12e-05} & {\tt   4.55}\\
                  $160 \times 160$ & {\tt   9.73e-04} & {\tt   1.93}& {\tt   1.68e-05} & {\tt   3.80}& {\tt   1.64e-06} & {\tt   5.26}& {\tt   3.99e-07} & {\tt   6.29}\\
               \end{tabular}\\[2ex]
         \end{center}
      \label{Tab: Shallow results final}
      \end{table}

      \begin{table}[ht]
       \centering
           \caption{\corr{Relative losses for the shallow-water equations, for varying tolerances reltol in the time-integration method {\tt lsode}, at the time when the solution becomes discontinuous.
                    Fourth-order discretization on a $20\times 20$ grid.}{Reviewer 2, Remarks 2 and 3: conservation for all examples}}
       \begin{center}
           \begin{tabular}{r|rrr|rrrrrrrrrr}
                  \multicolumn{1}{c}{} & \multicolumn{3}{c}{\bf Uniform grid} & \multicolumn{3}{c}{\bf Curvilinear grid}\\
             {\bf reltol} &  {\bf mass loss} & {\bf mom. loss} & {\bf energy loss} &  {\bf mass loss} & {\bf mom. loss} & {\bf energy loss} \\ \hline
              {\tt 1e-05} &   {\tt 8.07e-16} &      {\tt 1.49e-15} &    {\tt 1.16e-05} &   {\tt 2.02e-16} &      {\tt 5.40e-14} &    {\tt 3.35e-06} \\
              {\tt 1e-06} &   {\tt 1.41e-15} &      {\tt 9.39e-16} &    {\tt 8.07e-07} &   {\tt 4.04e-16} &      {\tt 5.55e-14} &    {\tt 4.85e-08} \\
              {\tt 1e-07} &   {\tt 0.00e+00} &      {\tt 7.04e-16} &    {\tt 2.07e-07} &   {\tt 8.07e-16} &      {\tt 5.56e-14} &    {\tt 1.59e-09} \\
              {\tt 1e-08} &   {\tt 0.00e+00} &      {\tt 1.37e-15} &    {\tt 5.05e-09} &   {\tt 2.02e-16} &      {\tt 5.55e-14} &    {\tt 7.83e-10} \\
              {\tt 1e-09} &   {\tt 1.01e-15} &      {\tt 5.25e-16} &    {\tt 1.36e-10} &   {\tt 4.04e-16} &      {\tt 5.44e-14} &    {\tt 1.05e-11} \\
              {\tt 1e-10} &   {\tt 0.00e+00} &      {\tt 1.11e-15} &    {\tt 3.91e-11} &   {\tt 0.00e+00} &      {\tt 5.68e-14} &    {\tt 3.12e-12} \\
              {\tt 1e-11} &   {\tt 6.06e-16} &      {\tt 9.39e-16} &    {\tt 2.12e-12} &   {\tt 1.01e-15} &      {\tt 5.48e-14} &    {\tt 5.79e-13} \\
              {\tt 1e-12} &   {\tt 4.04e-16} &      {\tt 8.30e-16} &    {\tt 5.61e-13} &   {\tt 6.06e-16} &      {\tt 5.53e-14} &    {\tt 4.02e-14} \\
         \end{tabular}\label{tab:conservation_shallow}
         \end{center}
      \end{table}

\section{Conclusion}
\label{sec:conclusion}
   This paper describes the construction of an arbitrary-order symmetry-preserving finite-difference technique on structured curvilinear staggered grids, offering flexibility and accuracy of the numerical approximations. The numerical examples presented in this work show that the method leads to results in which a high accuracy can be obtained, while the discrete mass, momentum and energy are all preserved. Also, the computation times are investigated for the linear-wave model, and the method is shown to be faster than when using a symmetry-preserving Galerkin-type approach. We have shown that the finite-difference method is able to handle difficult operators such as the advection operator in the shallow-water equations.

   This paper does not address the application of energy-conserving (symplectic) time integration. This is left outside the scope of this paper, and so are the discussion of boundary conditions and the use of numerical viscosity (for instance using flux limiters), which is necessary when shocks occur in the solution.  Finally, future work includes handling local grid refinements.
   
   {\bf Acknowledgments} The authors gratefully wish to acknowledge the useful comments provided by Jennifer Ryan (University of East Anglia) and Kees Vuik (Delft University of Technology) that helped to shape this work. 
We thank the reviewers for the suggestions they provided for revising our paper. 
   We also thank Shell and VORtech that allowed us to develop the research.

\bibliographystyle{plain}
\bibliography{References}

\bigskip
\Large{\textbf{Code Availability}}

\normalsize
\fbox{\begin{minipage}{0.94\textwidth}
The source code used for the experiments presented in this article can be obtained from:
\begin{center}
 \url{https://gitlab.com/VORtechBV/mamec}
\end{center}
\end{minipage}}
\end{document}